\documentclass[12pt]{article}
\usepackage{amssymb, amsfonts, amssymb, amsthm, amscd}

\newtheorem{theorem}{{\bf Theorem}}[section]
\newtheorem{lemma}{{\bf Lemma}}[section]
\newtheorem{proposition}{{\bf Proposition}}[section]
\newtheorem{definition}{{\bf Definition}}[section]
\newtheorem{corollary}{{\bf Corollary}}[section]
\newtheorem{remark}{{\bf Remark}}[section]
\newtheorem{example}{{\bf Example}}[section]

\begin{document}

\title{Reflected BSDE with a Constraint and its applications in incomplete market\footnote{This work is supported by the National
Basic Research Program of China (973 Program), No. 2007CB814902 and
No. 2007CB814906.}}
\author{Shige PENG$^{a,c}$\thanks{%
The author would like to acknowledge the partial support from the
Natural Science Foundation of China, grant No. 10131040. },\ \ \
Mingyu XU $^{b,c}$\footnote{Corresponding author, Email:
xvmingyu@gmail.com}
\\
{\small $^a$School of Mathematics and System Science, Shandong University, }%
\\
{\small 250100, Jinan, China}
\\
 {\small $^b$ Institute of Applied Mathematics, Academy of Mathematics and
Systems
Science,}\\
{\small Chinese Academy of Sciences, Beijing, 100080,
China.}\\{\small $^c$Department of Financial Mathematics and Control
science, School
of Mathematical Science,} \\
{\small Fudan University, Shanghai, 200433, China.}}
\date{The first version was submitted on arxiv on November, 2006.\\ This is the second version.}
\maketitle

\textbf{Abstract.} In this paper, we study a type of reflected BSDE
with a constraint and introduce a new kind of nonlinear expectation
via BSDE with a constraint and prove the Doob-Meyer decomposition
with respect to the super(sub)martingale introduced by this
nonlinear expectation. Then we an application on the pricing of
American options in incomplete market. \\[.5cm]

\section{Introduction}

El Karoui, Kapoudjian, Pardoux, Peng and Quenez (1997) studied the
problem of BSDE (backward stochastic differential equation) with
reflecting barrier, which is, a standard BSDE with an additional
continuous, increasing process in this equation to keep the solution
above a certain given continuous boundary process. This increasing
process must be chosen in certain minimal way, i.e. an integral
condition, called Skorokhod reflecting condition (cf. \cite
{Skoro1965}), is satisfied. The advantage of introducing the above
Skorokhod condition is that it possesses a very interesting coercive
structure which permits us to obtain many useful properties such as
uniqueness, continuous dependence and other kind of regularities. It
turns out to be a powerful tool to obtain the regularity properties
of the corresponding solutions of PDE with obstacle such as free
boundary PDE. Recently, this Skorokhod condition is generalized to
the case where the barrier $L$ is an $\mathbf{L}^{2}$--process in
\cite {PX2005}.

An important application of the constrained BSDE is the pricing of
contingent claims in an incomplete market, where the portfolios of an asset
is constrained in a given subset. In this case the solution $(y,z)$ of the
corresponding reflected BSDE must remain in this subset. In the pricing of
American options in the incomplete market, the related BSDE is a reflected
BSDE with constrained portfolios. This problem was studied by Karaztas and
Kou (cf. \cite{KK}). They required a condition that the constraint should be
a convex subset, the coefficient of the corresponding BSDE was also assumed
to be a linear, or at least a concave function. This limitation is mainly
due to the duality method applied as a main approach in that paper.

The main conditions of our paper is: $g$ is a Lipschitz function and the
constraint $\Gamma _{t}(\omega )$, $t\in [0,T]$ is a non--empty closed set.
The existence of such smallest $\Gamma $--constrained supersolution of BSDE
with coefficient $g$ is obtained in \cite{P99}. An interesting point of view
is that this supersolution is, in fact, the solution of the BSDE with a
singular coefficient $g_{\Gamma }$ defined by
\[
g_{\Gamma }(t,y,z)=g(t,y,z)1_{\Gamma _{t}}(y,z)+(+\infty )\cdot 1_{\Gamma
_{t}^{C}}(y,z).
\]
(see Remark \ref{singular} in appendix for details). One main result
of this paper, is the existence and uniqueness of reflected BSDEs
with this singular coefficient $g_{\Gamma }$ and we provide the
related generalized Skorokhod reflecting condition. Since our
coefficient $g$ as well as our constraint $\Gamma $ need not to be
concave or convex, the results of our paper provide a wide space of
freedom to treat different types of situations. Typically, in the
situation of differential games, the coefficients is neither convex
nor concave (see \cite{HL1995}, \cite{HL1995a} and \cite{HLP}).

Recent developments of continuous time finance requires a nonlinear
version of time consistent expectation. In 1997, the first author
has introduced a Brownian filtration $(\mathcal{F}_{t})_{t\geq 0}$
consistent nonlinear expectation
\[
\mathcal{E}^{g}[X]:X\in L^{2}(\Omega ,\mathcal{F}_{T},P)\rightarrow
\mathbb{R}
\]
call $g$--expectation, which is defined by $y_{0}^{X}$, where $%
(y_{t}^{X},z_{t}^{X})_{0\leq t\leq T}$ is the solution of the BSDE with a
given coefficient $g(t,y,z)$ and terminal condition $%
X $. Here we assume $g$ satisfies Lipschitz condition in $(y,z)$ as
well as $g(t,y,0)\equiv 0$. When $g$ is a linear function in
$(y,z)$, this $g$-expectation $\mathcal{E}^{g}[\cdot ]$ is just a
Girsanov transformation. But it becomes a nonlinear functional once
$g$ is nonlinear in $(y,z)$, i.e., $\mathcal{E}^{g}[\cdot ]$ is a
constant
preserving monotonic and nonlinear functional defined on $L^{2}(\Omega ,%
\mathcal{F}_{T},P)$.

Recently a profound link between super--replication, risk measures
(cf. \cite {ADEH1999}, \cite{Fo-Sc}) nonlinear expectations have
being explored (cf. \cite{B-El}, \cite{rosazza}, \cite{Peng03}). We
hope that the results of this paper will be proved to be useful in
this direction. We also refer to \cite{Dykin-Yushkevic},
\cite{Bensou1974}, \cite{Neveu1975}, \cite {Bismut1977}, \cite{E},
\cite{ALM1982}, \cite{Morimoto1984}, \cite{HL00} for interesting
research works in this domain.

To do researches for incomplete financial market, similarly as the
above $g$--expectation, we can also define the
corresponding $g_\Gamma $--expectation the smallest solution of BSDE with $%
g_\Gamma $ as well the corresponding $g_\Gamma $--supermartingales
and submartingales. We shall prove a $g_\Gamma $--supermartingale
decomposition theorem, which is a nonlinear version of Doob--Meyer
decomposition theorem. We point out that for the $g_\Gamma
$--submartingale decomposition can not be obtained by the above
mentioned $g_\Gamma $--supermartingale decomposition. We shall
obtain this decomposition theory in a quite different way.

This paper is organized as follows. In the next section we list our
main notations and main conditions required. In Section 3 we present the definition and some properties of $%
g_{\Gamma }$-expectation, with applications. In section 4, we prove
the results and proofs of the existence and uniqueness of reflected
BSDE with constraints. After introducing the definitions of $g_{\Gamma }$-martingale and $g_{\Gamma }$%
-super(sub)martingale, we prove the nonlinear Doob-Meyer's type
decomposition theorem corresponding to $g_{\Gamma
}$-super(sub)martingale in section 5. Then we give an application of
reflected BSDE with constraints: pricing of American option in
incomplete market in section 6. At last some useful results are
presented in appendix.

\section{$g_\Gamma $--solution: the smallest $g$-supersolution of BSDE
with constraint $\Gamma $}

Let $(\Omega ,\mathcal{F},P)$ be a probability space, and
$B=(B^1,B^2,\cdots ,B^d)^T$ be a $d$-dimensional Brownian motion
defined on $[0,\infty )$. We denote $\{\mathcal{F}_t;0\leq t<\infty
\}$ to be the natural filtration
generated by this Brownian motion $B:$%
\[
\mathcal{F}_t=\sigma \{\{B_s;0\leq s\leq t\}\cup \mathcal{N\}},
\]
where $\mathcal{N}$ is the collection of all $P-$null sets of
$\mathcal{F}$. The Euclidean norm of an element $x\in \mathbb{R}^m$
is denoted by $|x|$. We also need the following notations, for $p\in
[1,\infty )$:

\begin{itemize}
\item  $\mathbf{L}^{p}(\mathcal{F}_{t};\mathbb{R}^{m}):=$\{$\mathbb{R}^{m}$%
-valued $\mathcal{F}_{t}$--measurable random variables $X$ s.t. $%
E[|X|^{p}]<\infty $\};

\item  $\mathbf{L}_{\mathcal{F}}^{p}(0,t;\mathbb{R}^{m}):=$\{$\mathbb{R}^{m}$%
--valued and $\mathcal{F}_{t}$--progressively measurable processes $\varphi $ defined on $%
[0,t]$, s.t. $E\int_{0}^{t}|\varphi _{s}|^{p}ds<\infty $\};

\item  $\mathbf{D}_{\mathcal{F}}^{p}(0,t;\mathbb{R}^{m}):=$\{$\mathbb{R}^{m}$%
--valued and RCLL $\mathcal{F}_{t}$--progressively measurable processes $%
\varphi $ defined on $[0,t]$, s.t. $E[\sup_{0\leq s\leq t}|\varphi
_{s}|^{p}]<\infty $\};

\item  $\mathbf{A}_{\mathcal{F}}^{p}(0,t):=$\{increasing processes $A$ in $%
\mathbf{D}_{\mathcal{F}}^{p}(0,t;\mathbb{R})$ with $A(0)=0$\}.
\end{itemize}

\noindent When $m=1$, they are simplified as $\mathbf{L}^p(\mathcal{F}_t)$, $%
\mathbf{L}_{\mathcal{F}}^p(0,t)$ and
$\mathbf{D}_{\mathcal{F}}^p(0,t)$, respectively. We mainly interest
 the case of $p=2$. In this section, we consider BSDE on the interval
$[0,T]$, with a fixed $T>0$.

We consider a function
\[
g(\omega ,t,y,z):\Omega \times [0,T]\times \mathbb{R\times R}^{d}\rightarrow %
\mathbb{R}
\]
which always plays the role of the coefficient of our BSDE. $g$ satisfies
the following assumption: there exists a constant $\mu >0$, such that, for
each $y,y^{\prime }\;$in $\mathbb{R}$ and $z,z^{\prime }$ in $\mathbb{R}^{d}$%
, we have
\begin{equation}
\begin{array}{rll}
\mbox{(i)} & g(\cdot ,y,z)\in \mathbf{L}_{\mathcal{F}}^{2}(0,T); &  \\
\mbox{(ii)} & \left| g(t,\omega ,y,z)-g(t,\omega ,y^{\prime
},z^{\prime })\right| \leq \mu (\left| y-y^{\prime }\right| +\left|
z-z^{\prime }\right| ), & dP\times dt\mbox{ a.s. }
\end{array}
\label{Lip}
\end{equation}

Our constraint is described by $\Gamma (t,\omega ):\Omega \times
[0,T]\rightarrow
\mathcal{C}(\mathbb{R\times R}^{d})$, where $\mathcal{C}(\mathbb{R\times R}%
^{d})$ is the collection of all closed non--empty subsets of $%
\mathbb{R\times R}^{d}$, $\Gamma (t,\omega )$, which is
$\mathcal{F}_{t}$--adapted, namely,
\begin{equation}
\begin{array}{ll}
\mbox{(i)} & (y,z)\in \Gamma (t,\omega )\;\;\mbox{iff}\;\;d_{\Gamma
(t,\omega )}(y,z)=0,\;t\in [0,T]\mbox{, a.s.;} \\
\mbox{(ii)} & d_{\Gamma (\cdot )}(y,z)\mbox{ is }\mathcal{F}_{t}\mbox{%
--adapted process,\ for each }(y,z)\in \mathbb{R}\times
\mathbb{R}^{d},
\end{array}
\label{Gamma}
\end{equation}
where $d_{\Gamma _{\cdot }}(\cdot ,\cdot )$ is a distant function
from $(y,z)$ to $\Gamma $: for $t\in [0,T]$,
\[
d_{\Gamma _{t}}(y,z):=\inf_{(y^{\prime },z^{\prime })\in \Gamma
_{t}}(|y-y^{\prime }|^{2}+|z-z^{\prime }|^{2})^{1/2}\wedge 1.
\]
$d_{\Gamma _{t}}(y,z)$ is a Lipschitz function: for each $y,y^{\prime }\;$in
$\mathbb{R}$ and $z,z^{\prime }$ in $\mathbb{R}^{d}$, we always have
\[
\left| d_{\Gamma _{t}}(y,z)-d_{\Gamma _{t}}(y^{\prime },z^{\prime })\right|
\leq (|y-y^{\prime }|^{2}+|z-z^{\prime }|^{2})^{1/2}.
\]

\begin{remark}
The constraint discussed in \cite{P99} is
\begin{equation}
\Gamma _{t}(\omega )=\{(y,z)\in \mathbb{R}^{1+d}:\Phi (\omega ,t,y,z)=0\}.
\label{Fg}
\end{equation}
Here $\Phi (\omega ,t,y,z):\Omega \times [0,T]\times \mathbb{R}\times %
\mathbb{R}^{d}\rightarrow \mathbf{[}0,\infty )$ is a given nonnegative
measurable function, and satisfies integrability condition and Lipschitz
condition. In this paper we always consider the case
\[
\Phi (t,y,z)=d_{\Gamma _{t}}(y,z).
\]
\end{remark}

We are then within the framework of super(sub)solution of BSDE of the
following type:

\begin{definition}
($g$\textbf{--super(sub)solution}, cf. El Karoui, Peng and Quenez (1997)
\cite{EPQ} and Peng (1999) \cite{P99}) A process $y\in \mathbf{D}_{\mathcal{F%
}}^{2}(0,T)$ is called a $g$--supersolution (resp. $g$--subsolution), if
there exist a predictable process $z\in \mathbf{L}_{\mathcal{F}}^{2}(0,T;%
\mathbb{R}^{d})$ and an increasing RCLL process $A\in \mathbf{A}_{\mathcal{F}%
}^{2}(0,T)$ (resp. $K\in \mathbf{A}_{\mathcal{F}}^{2}(0,T)$), such that $%
t\in [0,T],$%
\begin{eqnarray}
y_{t}
&=&y_{T}+\int_{t}^{T}g(s,y_{s},z_{s})ds+A_{T}-A_{t}-\int_{t}^{T}z_{s}dB_{s},
\label{CBSDE} \\
\mbox{(resp. }y_{t}
&=&y_{T}+\int_{t}^{T}g(s,y_{s},z_{s})ds-(K_{T}-K_{t})-%
\int_{t}^{T}z_{s}dB_{s}.\mbox{)}  \nonumber
\end{eqnarray}
Here $z$ and $A$ (resp. $K$) are called the martingale part and increasing
part, respectively. $y$ is called a $g$--solution if $A_{t}=K_{t}=0$, for $%
t\in [0,T]$. $y$ is called a $\Gamma $--constrained $g$--supersolution if $y$
and its corresponding martingale part $z$ satisfy
\begin{equation}
(y_{t},z_{t})\in \Gamma _{t},\,\,\,(\hbox{or }d_{\Gamma
_{t}}(y_{t},z_{t})=0),\,\,\,\,dP\times dt\hbox{ a.s. in }\Omega \times [0,T],
\label{G-constr}
\end{equation}
\end{definition}

\begin{remark}
We observe that, if $y\in \mathbf{D}_{\mathcal{F}}^{2}(0,T)$ is a $g$%
--supersolution or $g$--subsolution, then the pair $(z,A)$ in (\ref{CBSDE})
are uniquely determined since the martingale part $z$ is uniquely
determined.\ Occasionally, we also call $(y,z,A)$ a $g$--supersolution or $g$%
--subsolution.
\end{remark}

By \cite{P99}, (see Appendix Theorem \ref{mono1}), if there exists at least
one $\Gamma $--constrained $g$--supersolution, then the smallest $\Gamma $%
--constrained $g$--supersolution exists. In fact, a $\Gamma $-constraint $g$%
-supersolution can be considered as a solution of the BSDE with a singular
coefficient $g_{\Gamma }$ defined by
\[
g_{\Gamma }(t,y,z)=g(t,y,z)1_{\Gamma _{t}}(y,z)+(+\infty )\cdot 1_{\Gamma
_{t}^{C}}(y,z).
\]
So we define the smallest $\Gamma $--constrained $g$--supersolution by $%
g_{\Gamma }$\textbf{--}solution.

\begin{definition}
\label{gGsol}\textbf{(}$g_{\Gamma }$\textbf{--solution)} $y$ or $%
(y_{t},z_{t},A_{t})_{0\leq t\leq T}$ is called $g_{\Gamma
}$--solution on $[0,T]$ with a given terminal condition $X$ if it is
the smallest $\Gamma $--constrained $g$--supersolution with
$y_{T}=X$:
\begin{eqnarray}
y_{t}
&=&X+\int_{t}^{T}g(s,y_{s},z_{s})ds+A_{T}-A_{t}-\int_{t}^{T}z_{s}dB_{s},
\label{gG-solution} \\
\ d_{\Gamma _{t}}(y_{t},z_{t}) &=&0,\,\,\,\,dP\times dt\hbox{ a.s. in }%
\Omega \times [0,T],\;\;dA_{t}\geq 0,\;\;t\in [0,T].  \nonumber
\end{eqnarray}
In other words, if there exists another triple $(y',z',A')$
satisfying (\ref{gG-solution}), then $y'_t \geq y_t$, for
$t\in[0,T]$.
\end{definition}

\begin{remark}
\label{Dsmal} The above definition does not imply that the
increasing process $A$ is also the smallest one, i.e. for another
triple $(\bar{y},\bar{z},\bar{A})$ satisfying (\ref{gG-solution}),
we may have $A_{t}\geq \bar{A}_{t}$.
\end{remark}

An example is as following.

\begin{example}
Consider the case when $[0,T]=[0,2]$, $X=0$, $g=0$ and $\Gamma
_{t}=\{(y,z):y\geq 1_{[0,1]}(t)\}$. So the $g_{\Gamma }$-solution of
this equation is the solution of reflected BSDE with lower barrier
$1_{[0,1]}(t)$. It's easy to
see that the smallest solution is $y_{t}=1_{[0,1)}(t)$ with $z_{t}=0$, $%
A_{t}=1_{[1,2]}(t)$. Obviously $\overline{y}_{t}=1_{[0,2)}(t)$ with $%
\overline{z}_{t}=0$, $\overline{A}_{t}=1_{\{t=2\}}(t)$ is another $\Gamma $%
--constrained $g$--supersolution with the same terminal condition $%
y_{T}^{\prime }=0$. However we have $A_{t}>\overline{A}_{t}$ on the interval
$[1,2)$.
\end{example}

\section{Nonlinear Expectation: $g_{\Gamma }$-expectation and its properties}

In this section we first introduce a new type of
$\mathcal{F}$--consistent nonlinear expectations via $g_{\Gamma
}$--solutions, then we study the properties of this nonlinear
expectations. At last an application for risk measure in the incomplete market
 is concerned. We assume: there exists a large enough constant $C_{0}$ such that for $%
\forall y\geq C_{0}$%
\begin{equation}
g(t,y,0)\leq C_{0}+\mu |y|,\;\mbox{and \ }(y,0)\in \Gamma _{t};
\label{g-Gamma}
\end{equation}
and the terminal conditions to be in the following linear subspace of $%
\mathbf{L}^{2}(\mathcal{F}_{T})$:
\[
\mathbf{L}_{+,\infty }^{2}(\mathcal{F}_{T}):=\{\xi \in \mathbf{L}^{2}(%
\mathcal{F}_{T}),\;\xi ^{+}\in \mathbf{L}^{\infty
}(\mathcal{F}_{T})\}.
\]

\begin{proposition}
\label{smal}We assume (\ref{Lip}), (\ref{Gamma}) and (\ref{g-Gamma}) hold. Then for each $%
X\in \mathbf{L}_{+,\infty }^{2}(\mathcal{F}_{T})$, the $g_{\Gamma
}$-solution with terminal condition $y_{T}=X$ exists. Furthermore,
we have $y_{t}\in \mathbf{L}_{+,\infty }^{2}(\mathcal{F}_{t})$, for
$t\in [0,T]$.
\end{proposition}

\noindent\textbf{Proof. }We consider
\[
y_{0}(t)=(\left\| X^{+}\right\| _{\infty }\vee C_{0})e^{\mu
(T-t)}+C_{0}(T-t)+(X-\left\| X^{+}\right\| _{\infty }\vee
C_{0})1_{\{t=T\}}.
\]
It is the solution of the following backward equation:
\[
y_{0}(t)=X+\int_{t}^{T}(C_{0}+\mu |y_{0}(s)|)ds+A^{0}(T)-A^{0}(t),
\]
where $A^{0}$ is an increasing process: $A^{0}(t):=(\left\|
X^{+}\right\| _{\infty }\vee C_{0}-X)1_{t=T}$. Meanwhile
$y_{0}(\cdot )$ can be expressed as:
\[
y_{0}(t)=X+\int_{t}^{T}g(s,y_{0}(s),0)ds+\int_{t}^{T}[c+\mu
|y_{0}(s)|-g(s,y_{0}(s),0)]ds+A^{0}(T)-A^{0}(t).
\]
Thus the triple defined on $[0,T]$ by
\[
(y_{1}(t),z_{1}(t),A_{1}(t)):=(y_{0}(t),0,\int_{0}^{t}[c+\mu
|y_{0}(s)|-g(s,y_{0}(s),0)]ds+A^{0}(t))
\]
is a $\Gamma $--constrained $g$--supersolution with $y_{1}(T)=X$.
According to Theorem \ref{exist} in appendix, the $g_{\Gamma
}$--solution with $y(T)=X$ exists. We
also have $(y_{t})^{+}\in \mathbf{L}^{\infty }(\mathcal{F}_{T})$ since $%
y_{t}\leq y_{1}(t)=y_{0}(t)$. \hfill $\Box $\medskip

We now introduce the notion of $g_{\Gamma }$--expectation:

\begin{definition}
\label{eval}We assume that for each $0\leq t\leq T<\infty $,
$g(t,0,0)=0$ and $(0,0)\in \Gamma _{t}$, assumptions (\ref{Lip}),
(\ref{Gamma}) and (\ref
{g-Gamma})hold. Then consider $X\in \mathbf{L}_{+,\infty }^{2}(\mathcal{F}_{T})$%
, let $(y,z,A)$ be the $g_{\Gamma }$--solution defined on $[0,T]$
with terminal condition $y_{T}=X$. We define
$\mathcal{E}_{t,T}^{g_{\Gamma }}[X]:=y_{t}.$ The system
\[
\mathcal{E}_{t,T}^{g_{\Gamma }}[\cdot ]:\mathbf{L}_{+,\infty }^{2}(\mathcal{F%
}_{T})\rightarrow \mathbf{L}_{+,\infty }^{2}(\mathcal{F}_{t}),\,\,\,\,\,\,0%
\leq t\leq T<\infty
\]
is called $g_{\Gamma }$-expectation.
\end{definition}

\begin{remark}
Under assumptions (\ref{Lip}), (\ref{Gamma}) and (\ref {g-Gamma}),
proposition \ref{smal} guarantees the existence of $g_{\Gamma
}$-expectation.
\end{remark}

We have

\begin{proposition}
\label{eval1}A $g_{\Gamma }$-expectation is an
$\mathcal{F}$-consistent expectation, i.e., it satisfies the
followings: for each $0\leq t\leq
T<\infty $ and $X$, $X^{\prime }\in \mathbf{L}_{+,\infty }^{2}(\mathcal{F}%
_{T})$,

\textbf{(A1) }\textsl{Monotonic property}\textbf{:\ }$\mathcal{E}%
_{t,T}^{g_{\Gamma }}[X]\leq \mathcal{E}_{t,T}^{g_{\Gamma
}}[X^{\prime }],\,\,\,\,$ if $X\leq X^{\prime }$;\

\textbf{(A2) }\textsl{Self-preserving:}\textbf{\ }$\mathcal{E}%
_{T,T}^{g_{\Gamma }}[X]=X$;

\textbf{(A3) }\textsl{Time consistency:}\textbf{\ }$\mathcal{E}%
_{s,t}^{g_{\Gamma }}[\mathcal{E}_{t,T}^{g_{\Gamma }}[X]]=\mathcal{E}%
_{s,T}^{g_{\Gamma }}[X]$, $\,\,\,\,\,0\leq s\leq t\leq T$;

\textbf{(A4) }\textsl{1-0 law: }$1_{D}\mathcal{E}_{t,T}^{g_{\Gamma }}[X]=%
\mathcal{E}_{t,T}^{g_{\Gamma }}[1_{D}X]$, $\,\,\,\,\forall D\in \mathcal{F}%
_{t}$.
\end{proposition}

\noindent\textbf{Proof. } \textbf{(A1)} is a direct consequence of
the comparison theorem \ref{comp} of the $g_{\Gamma }$--solution.
\textbf{(A2) }is obvious. For \textbf{(A3)}, it is easy the
check that, if $(y_{s})_{0\leq s\leq T}$ is the $g_{\Gamma }$--solution on $%
[0,T]$ with $y_{T}=X$, then $(y_{s})_{0\leq s\leq t}$ is also the
$g_{\Gamma }$--solution on $[0,t]$ with the fixed terminal condition
$y_{t}$.

To prove \textbf{(A4)}, we multiply $1_{D}$ to two sides of the
equation, for $t\leq s\leq T$, since $g(s,0,0)\equiv 0$, and
$d_{\Gamma _{s}}(0,0)\equiv 0$, we have
\begin{eqnarray*}
1_{D}y_{s}
&=&1_{D}X+\int_{s}^{T}g(r,1_{D}y_{r},1_{D}z_{r})dr+1_{D}A_{T}-1_{D}A_{s}-%
\int_{s}^{T}1_{D}z_{r}dB_{r}, \\
&&d_{\Gamma _{s}}(1_{D}y_{s},1_{D}z_{s})\equiv 0.
\end{eqnarray*}
Thus it is obvious that $(1_{D}y_{s},1_{D}z_{s})_{t\leq s\leq T}$
must be the $g_{\Gamma }$--solution on $[s,T]$ with $y_{T}1_{D}$ as
the terminal condition, which implies \textbf{(A4)}.\hfill $\Box
$\medskip

Moreover, by the comparison theorem for $g_{\Gamma }$-solution, we
have

\begin{proposition}
Under assumptions (\ref{Lip}), (\ref{Gamma}) and (\ref{g-Gamma}), for each $%
0\leq t\leq T<\infty $ and $X\in \mathbf{L}_{+,\infty }^{2}(\mathcal{F}_{T})$%
,\textbf{\ }if
\[
\Gamma _{t}^{1}\supseteq \Gamma _{t}^{2}\,\,\,\mbox{and\ \ \ }%
g^{1}(t,y,z)\leq g^{2}(t,y,z),\,
\]
then $\mathcal{E}_{t,T}^{g_{\Gamma ^{1}}^{1}}[X]\leq \mathcal{E}%
_{t,T}^{g_{\Gamma ^{2}}^{2}}[X]$.
\end{proposition}

Now we study some properties of $g_{\Gamma }$-expectation associated
with dynamic risk measure, such as constant preserving property,
positive homogenous property, convex property, sublinear property,
constant translation invariant property and subadditive property.
And in the following of this section, we always assume that
assumptions (\ref{Lip}) and (\ref{Gamma}) hold.

\begin{proposition}[positive homogenous and convexity]
\label{cons}If $g(t,y,0)=0$, and $\mathbb{R\times \{}0\}\subset \Gamma _{t}$%
, $t\in [0,T]$, then $g_{\Gamma }$-expectation is conditional
constant preserving,
\[
\mathcal{E}_{t,T}^{g_{\Gamma }}[X]=X,\mbox{ for }X\in
\mathbf{L}_{+,\infty }^{2}(\mathcal{F}_{t}).
\]
Specially, for $C\in \mathbb{R}$, $\mathcal{E}_{t,T}^{g_{\Gamma
}}[C]=C.$
\end{proposition}

\noindent\textbf{Proof. } For $X\in \mathbf{L}_{+,\infty
}^{2}(\mathcal{F}_{t})$, it is easy to check that
$(y_{t},z_{t},A_{t})\equiv (X,0,0)$ is the $g_{\Gamma }$-solution of
constraint BSDE associated to $(X,g,\Gamma )$, in view of $g(t,y,0)=0$, and $%
\mathbb{R\times \{}0\}\subset \Gamma _{t}$, $t\in [0,T]$. So the
result follows. \hfill$\square $

\begin{proposition}
\label{HoCo}Set $g(t,0,0)=0$ and $(0,0)\in \Gamma _{t}$ hold for
each $0\leq t\leq T<\infty $,

(i) under assumption (\ref{g-Gamma}), the nonlinear
$\mathcal{F}$-consistent expectation, $g_{\Gamma }$-expectation is
positive homogenous, i.e.
\[
\mathcal{E}_{t,T}^{g_{\Gamma }}[cX]=c\mathcal{E}_{t,T}^{g_{\Gamma }}[X],%
\mbox{ for }c>0,X\in \mathbf{L}_{+,\infty }^{2}(\mathcal{F}_{T}),
\]
if $g$ is positive homogenous in $(y,z)$ and $\Gamma _{t}$ is a cone for $t\in [0,T]$, i.e. if $%
(y,z)\in \Gamma _{t}$, then for $c>0$, $(cy,cz)\in \Gamma _{t}$;

(ii) under assumption (\ref{g-Gamma}), if $g$ and $\Gamma $ are
convex in $(y,z)$, then $g_{\Gamma }$-expectation is convex,
\[
\mathcal{E}_{t,T}^{g_{\Gamma }}[\alpha X_{1}+(1-\alpha )X_{2}]\leq
\alpha
\mathcal{E}_{t,T}^{g_{\Gamma }}[X_{1}]+(1-\alpha )\mathcal{E}%
_{t,T}^{g_{\Gamma }}[X_{2}],\mbox{ for }\alpha \in
[0,1],\;X_{1},X_{2}\in \mathbf{L}_{+,\infty }^{2}(\mathcal{F}_{T}).
\]
\end{proposition}

\noindent\textbf{Proof. }
(i) It is easy to see that $cX\in \mathbf{L}_{+,\infty }^{2}(\mathcal{F}%
_{T}) $, with $c>0$, if and only if $X\in \mathbf{L}_{+,\infty }^{2}(%
\mathcal{F}_{T})$. Let $(y,z,A)$ be the $g_{\Gamma }$--solution defined on $%
[t,T]$ with terminal condition $y_{T}=X$, i.e. for $t\leq s\leq T$,
\begin{eqnarray*}
y_{s}
&=&X+\int_{s}^{T}g(r,y_{r},z_{r})dr+A_{T}-A_{s}-\int_{s}^{T}z_{r}dB_{r},\; \\
d_{\Gamma _{s}}(y_{s},z_{s}) &=&0,\;\;\mbox{a.s., a.e.}
\end{eqnarray*}
Since $g$ is homogenous and $\Gamma $ is a cone, we have, for $c>0$, $%
(cy_{s},cz_{s})\in \Gamma _{s}$, a.s.a.e. and
\begin{eqnarray*}
cy_{s}
&=&cX+c\int_{s}^{T}g(r,y_{r},z_{r})dr+cA_{T}-cA_{s}-c\int_{s}^{T}z_{r}dB_{r}
\\
&=&cX+\int_{s}^{T}g(r,cy_{r},cz_{r})dr+cA_{T}-cA_{s}-c%
\int_{s}^{T}z_{r}dB_{r},
\end{eqnarray*}
It is obvious that $(cy,cz,cA)$ is the $g_{\Gamma }$--solution with
terminal condition $cX$, i.e. $\mathcal{E}_{t,T}^{g_{\Gamma }}[cX]=cy_{t}=c%
\mathcal{E}_{t,T}^{g_{\Gamma }}[X]$.

(ii) Since $X_{1},X_{2}\in \mathbf{L}_{+,\infty
}^{2}(\mathcal{F}_{T})$, and
$\alpha \in [0,1]$, so $\alpha X_{1}+(1-\alpha )X_{2}\in \mathbf{L}%
_{+,\infty }^{2}(\mathcal{F}_{T})$. We denote
$\mathcal{E}_{t,T}^{g_{\Gamma
}}[\alpha X_{1}+(1-\alpha )X_{2}]=y_{t}$, which is the $g_{\Gamma }$%
-solution of BSDE$(g,\Gamma )$ on $[t,T]$, with terminal condition
$\alpha X_{1}+(1-\alpha )X_{2}$, i.e. for $t\leq s\leq T$,
\begin{eqnarray}
y_{s} &=&\alpha X_{1}+(1-\alpha
)X_{2}+\int_{s}^{T}g(r,y_{r},z_{r})dr+A_{T}-A_{s}-\int_{s}^{T}z_{r}dB_{r},
\label{conv1} \\
d_{\Gamma _{s}}(y_{s},z_{s}) &=&0,\;\;\mbox{a.s.a.e.}  \nonumber
\end{eqnarray}
Set $\mathcal{E}_{t,T}^{g_{\Gamma }}[X_{1}]=y_{t}^{1}$ and $\mathcal{E}%
_{t,T}^{g_{\Gamma }}[X_{2}]=y_{t}^{2}$, where for $i=1,2$, $%
(y^{i},z^{i},A^{i})$ is the $g_{\Gamma }$-solution of BSDE with
terminal value $X^{i}$, associated to $(g,\Gamma )$, i.e.
\[
y_{s}^{i}=X^{i}+\int_{s}^{T}g(r,y_{r}^{i},z_{r}^{i})dr+A_{T}^{i}-A_{t}^{i}-%
\int_{s}^{T}z_{r}^{i}dB_{r},\;d_{\Gamma _{s}}(y_{t}^{i},z_{t}^{i})=0,\;\;%
\mbox{a.s.a.e.}
\]
Then we know that the convex combination $(\alpha y^{1}+(1-\alpha
)y^{2},\alpha z^{1}+(1-\alpha )z^{2},\alpha A^{1}+(1-\alpha )A^{2})$ is a $g$%
-supersolution of BSDE with terminal value $\alpha X_{1}+(1-\alpha
)X_{2}$ and coefficient $\widetilde{g}$, where
\[
\widetilde{g}(s,y,z)=\alpha g(r,y_{s}^{1},z_{s}^{1})+(1-\alpha )g(s,\frac{1}{%
1-\alpha }(y-\alpha y_{s}^{1}),\frac{1}{1-\alpha }(z-\alpha
z_{s}^{1})).
\]
Moreover Since $\Gamma _{s}$ is convex for $s\in [t,T]$, $(\alpha
y_{s}^{1}+(1-\alpha )y_{s}^{2},\alpha z_{s}^{1}+(1-\alpha
)z_{s}^{2})\in \Gamma _{s}$, a.s. a.e.. Notice that $g$ is a convex
function, we have
\begin{eqnarray*}
\widetilde{g}(s,y_{s},z_{s}) &=&\alpha
g(s,y_{s}^{1},z_{s}^{1})+(1-\alpha
)g(s,\frac{1}{1-\alpha }(y_{s}-\alpha y_{s}^{1}),\frac{1}{1-\alpha }%
(z_{s}-\alpha z_{s}^{1})) \\
&\geq &g(s,y_{s},z_{s}).
\end{eqnarray*}
By comparison theorem, and remember that $y_{t}$ is the $g_{\Gamma }$%
-solution, then
\[
\mathcal{E}_{t,T}^{g_{\Gamma }}[\alpha X_{1}+(1-\alpha
)X_{2}]=y_{t}\leq \alpha y_{t}^{1}+(1-\alpha )y_{t}^{2}=\alpha
\mathcal{E}_{t,T}^{g_{\Gamma }}[X_{1}]+(1-\alpha
)\mathcal{E}_{t,T}^{g_{\Gamma }}[X_{2}].
\]
$\square $

\begin{corollary}\label{sublin}[Sublinear] Let $g(t,0,0)=0$ and $(0,0)\in \Gamma _{t}$ hold, for each
$0\leq t\leq T<\infty $. If $g$ is sublinear in $(y,z)$, i.e. $g$ is
homogenous and subadditive in $(y,z)$, which implies for $c>0$,
$(y,z)$ and $(y',z')$ in $\mathbb{R}^{1+d}$,
\[
g(t,cy,cz)=cg(t,y,z) \mbox{ and }g(t,y+y',z+z')\leq
g(t,y,z)+g(t,y'z'),
\]
and $\Gamma _{t}$ is a convex cone for $t\in [0,T]$, then $g_{\Gamma
}$-expectation is sublinear.
\end{corollary}

\noindent\textbf{Proof. } Since sublinearity is equivalent to
convexity plus positive homogeneity, the thesis follows from
Proposition \ref{HoCo}. $\square $

\begin{proposition}[constant translation invariant]
\label{tr-con}For each $0\leq t\leq T<\infty $,

(i) if $g$ and $\Gamma $ only depend on $z$, $g(t,z)$ is bounded and
$0\in \Gamma _{t}$, then $g_{\Gamma }$-expectation is translation
invariant,
\[
\mathcal{E}_{t,T}^{g_{\Gamma }}[X+\eta
]=\mathcal{E}_{t,T}^{g_{\Gamma
}}[X]+\eta ,\mbox{ for }\eta \in \mathbf{L}_{+,\infty }^{2}(\mathcal{F}%
_{t}),X\in \mathbf{L}_{+,\infty }^{2}(\mathcal{F}_{T});
\]

(ii) if $g(t,y,z)=g_{1}(t,z)+ay$ with $g_{1}(t,z)$ is bounded and
$\Gamma $ only depends on $z$, with $0\in \Gamma _{t}$, then
$g_{\Gamma }$-expectation is constant invariant with discount factor
$e^{a(T-t)}$,
\[
\mathcal{E}_{t,T}^{g_{\Gamma }}[X+\eta
]=\mathcal{E}_{t,T}^{g_{\Gamma
}}[X]+\eta e^{a(T-t)},\mbox{ for }\eta \in \mathbf{L}_{+,\infty }^{2}(%
\mathcal{F}_{t}),X\in \mathbf{L}_{+,\infty }^{2}(\mathcal{F}_{T}).
\]
\end{proposition}

\noindent\textbf{Proof. } Obviously (\ref{g-Gamma}) is satisfied
under the assumption (i) and (ii).

(i) Obviously $X+\eta \in
\mathbf{L}_{+,\infty }^{2}(\mathcal{F}_{T})$. By the definition of $%
g_{\Gamma }$-expectation, we know that $\mathcal{E}_{t,T}^{g_{\Gamma
}}[X]:=y_{t}$, where $(y,z,A)$ is the $g_{\Gamma }$--solution of
constraint BSDE$(X,g,\Gamma )$ on $[t,T]$. So for $s\in [t,T]$,
\begin{eqnarray*}
y_{s}+\eta &=&X+\eta
+\int_{s}^{T}g(r,z_{r})dr+A_{T}-A_{s}-\int_{s}^{T}z_{r}dB_{r},\; \\
d_{\Gamma _{s}}(z_{s}) &=&0,\;\;\mbox{a.s.a.e.}
\end{eqnarray*}
It follows that $\mathcal{E}_{t,T}^{g_{\Gamma }}[X+\eta ]=y_{t}+\eta =%
\mathcal{E}_{t,T}^{g_{\Gamma }}[X]+\eta $.

(ii) By the definition of $%
g_{\Gamma }$-expectation, we know that $\mathcal{E}_{t,T}^{g_{\Gamma
}}[X]:=y_{t}$, where $(y,z,A)$ is the $g_{\Gamma }$--solution on
$[t,T]$ with terminal condition $y_{T}=X$. Since $\int_{s}^{T}a\eta
e^{a(T-r)}dr=\eta \int_{s}^{T}d(-e^{a(T-r)})=-\eta +\eta
e^{a(T-s)}$, we get
\begin{eqnarray*}
y_{s}+\eta e^{a(T-s)} &=&X+\eta
+\int_{s}^{T}[g_{1}(r,z_{r})+a(y_{r}+\eta
e^{a(T-r)})]dr+A_{T}-A_{s}-\int_{s}^{T}z_{r}dB_{r} \\
&=&X+\eta +\int_{s}^{T}g(r,y_{r}+\eta
e^{a(T-r)},z_{r})dr+A_{T}-A_{s}-\int_{s}^{T}z_{r}dB_{r}.
\end{eqnarray*}
Notice that we still have $d_{\Gamma _{s}}(z_{s})=0,\;\;$a.s.a.e..
And it is
easy to check that $(y,z,A)$ is the $g_{\Gamma }$--solution. Then $%
y_{s}+\eta e^{a(T-s)}$ is the $g_{\Gamma }$-solution of constraint BSDE$%
(X+\eta ,g,\Gamma )$, i.e.
\[
\mathcal{E}_{t,T}^{g_{\Gamma }}[X+\eta ]=y_{t}+\eta e^{a(T-t)}=\mathcal{E}%
_{t,T}^{g_{\Gamma }}[X]+\eta e^{a(T-t)}.
\]
$\square $

As we know from Rosazza \cite{rosazza}, we can use $g$-expectation
to describe risk measure dynamically. However in incomplete market,
since portfolio is constraint, risk of a financial position must
increase. This indicates us to use our $g_{\Gamma }$-expectation to
study dynamic risk measure in incomplete market.

\begin{example}[Risk measure with no-shortselling constraint] Set $\Gamma$ only depends on $z$, with $\Gamma_t
=\mathbb{R}^d_+$, and $g$ is Lipschitz in $(y,z)$, then for a
financial position $X\in \mathbf{L}^2_{+,\infty}({\cal F}_T)$ define
a dynamic risk measure:
\[
\rho_t(X)=\mathcal{E}_{t ,T }^{g_{\Gamma }}[-X].
\]
Thanks to Proposition \ref{HoCo}, \ref{tr-con}, and Corollary
\ref{sublin}, we have
\begin{itemize}
  \item $\rho_t(\cdot)$ is a dynamic convex
time-consistent risk measure, if $g$ is convex in $(y,z)$.
  \item $\rho_t(\cdot)$ is a dynamic coherent
time-consistent risk measure, if $g$ only depends on $z$ and is
sublinear in $z$.
  \item $\rho_t(\cdot)$ is a dynamic sublinear
time-consistent risk measure, if $g$ is sublinear in $(y,z)$.
\end{itemize}
If we define another dynamic risk measure $\bar{\rho}_t$, for a
financial position $X\in \mathbf{L}^2_{+,\infty}({\cal F}_T)$, by
\[
\bar{\rho}_t(X)=\mathcal{E}_{g}[-X|\mathcal{F}_t].
\]
Here $\mathcal{E}_{g}[\cdot|\mathcal{F}_t]$ is a $g$-expectation,
(cf. \cite{Peng03}). By comparison theorem for BSDE, we can easily
get
\[\rho_t(X)\geq \bar{\rho}_t(X),
\]
which implies that in the market with no-shortselling constraint,
for same financial position, we need more money to cover its risk.
\end{example}

\section{$g_{\Gamma }$--reflected BSDEs}

 Before we go further to study more properties of $g_{\Gamma
 }$-expectation, we change our attentions to $g_{\Gamma }$--reflected
 BSDEs, which will play important roles in further research.

\subsection{Existence of $g_{\Gamma }$--reflected BSDEs}

In this section we consider the smallest $g$--supersolution with constraint $%
\Gamma $ and a lower (resp. upper) reflecting obstacle $L$ (resp. $U$). We
assume that the two reflected obstacles $L$ and $U\ $are $\mathcal{F}_{t}$%
-adapted processes satisfying
\begin{equation}
L,\;U\in \mathbf{L}_{\mathcal{F}}^{2}(0,T)\;\;\;\;\;\hbox{ and }%
\,\,\,\,\,ess\sup_{0\leq t\leq T}L_{t}^{+},\;ess\sup_{0\leq t\leq
T}U_{t}^{-}\in \mathbf{L}^{2}(\mathcal{F}_{T}).\;  \label{barc}
\end{equation}

Here we focus on the constraint $\Gamma $ which does not depend on
$y$, only depends on $z$, i.e. $\Gamma (t,\omega ):\Omega \times
[0,T]\rightarrow
\mathcal{C}(\mathbb{R}^{d})$, where $\mathcal{C}(\mathbb{R}%
^{d})$ is the collection of all closed non--empty subsets of $%
\mathbb{R}^{d}$ and $\Gamma (t,\omega )$ is
$\mathcal{F}_{t}$--adapted. In fact, this condition of $\Gamma $ is
not an essential difficulty in following proofs in this section. We
can easily generalize the results to the case when also depends on
$y$.

First let us introduce the definition of $g_{\Gamma }$-reflected
solutions:

\begin{definition}
\label{dif-RCBSDE2}A $g_{\Gamma }$-reflected solution with a lower
obstacle $L$ is a quadruple of processes $(y,z,A,\bar{A})$
satisfying \newline (\textbf{i}) $(y,z,A,\bar{A})\in
\mathbf{D}_{\mathcal{F}}^{2}(0,T)\times
\mathbf{L}_{\mathcal{F}}^{2}(0,T;\mathbb{R}^{d})\times (\mathbf{A}_{\mathcal{%
F}}^{2}(0,T))^{2}$ verifies
\begin{eqnarray}
y_{t} &=&X+\int_{t}^{T}g(s,y_{s},z_{s})ds+A_{T}-A_{t}+\overline{A}_{T}-%
\overline{A}_{t}-\int_{t}^{T}z_{s}dB_{s},  \label{RCBSDEl} \\
&&\ \ \
\begin{array}{c}
d_{\Gamma _{t}}(z_{t})=0\mbox{, \ }\;\;\;dP\times dt\;\mbox{a.s..}
\end{array}
\nonumber
\end{eqnarray}
(\textbf{ii}) $y_{t}\geq L_{t}$ and the generalized Skorokhod
reflecting condition is satisfied: for each $L^{*}\in
\mathbf{D}_{\mathcal{F}}^{2}(0,T)$ such that $y_{t}\geq
L_{t}^{*}\geq L_{t}$, $dP\times dt$-a.s., we have
\begin{equation}
\int_{0}^{T}(y_{s-}-L_{s-}^{*})d\overline{A}_{s}=0,\;\mbox{a.s.,\ \
} \label{Skoro1}
\end{equation}
\textbf{(iii) }$y$ is the smallest one, i.e., for any quadruple
$(y^{*},z^{*},A^{*},\bar{A}^{*})$ satisfying \textbf{(i) }and
\textbf{(ii)}, we have
\[
y_{t}\leq y_{t}^{*},\;\;\;\;\;\forall t\in [0,T],\;\;\mbox{a. s.}.
\]
\end{definition}

Here we use two increasing processes $A$, $\overline{A}$ to push $y$
in order to keep the solution $(y,z)$ staying in constraint $\Gamma
$ and upper the barrier $L$ respectively. More precisely, the role
of $A$ is to keep the process $z$ staying in the given constraint
$\Gamma $, while $\overline{A}$ acts only when $y$ tends to cross
downwards the barrier $L$.

Our first main result in this section is:

\begin{theorem}
\label{existrl}Suppose (\ref{Lip}), (\ref{Gamma}) and (\ref{barc})
hold. For a given terminal condition $X\in
\mathbf{L}^{2}(\mathcal{F}_{T})$, we assume
that there exists a triple $(y^{*},z^{*},A^{*})\in \mathbf{D}_{\mathcal{F}%
}^{2}(0,T)\times \mathbf{L}_{\mathcal{F}}^{2}(0,T)\times \mathbf{A}_{%
\mathcal{F}}^{2}(0,T)$, such that $dA^{*}\geq 0$ and following hold
\begin{eqnarray}
y_{t}^{*}
&=&X+\int_{t}^{T}g(s,y_{s}^{*},z_{s}^{*})ds+(A_{T}^{*}-A_{t}^{*})-%
\int_{t}^{T}z_{s}^{*}dB_{s},  \label{H1} \\
&&\ \ \ \
\begin{array}{c}
(y_{t}^{*},z_{t}^{*})\in \Gamma _{t}\cap \{[L_{t},\infty )\times %
\mathbb{R}^{d}\},\;\;dP\times dt\mbox{-a.s.}.
\end{array}
\nonumber
\end{eqnarray}
Then there exists the $g_{\Gamma }$--reflected solution
$(y,z,A,\bar{A})$ with the barrier $L$ of Definition
\ref{dif-RCBSDE2}.
\end{theorem}

\begin{remark}
This theorem can be generalized to the case when $\Gamma $ also depends on $y$%
easily.
\end{remark}

The smallest $g_{\Gamma }$--reflected solution with a upper obstacle
$U$ is relatively more complicated than the case of the lower
obstacle.

\begin{definition}
\label{RCBSDE1}The $g_{\Gamma }$--reflected solution with an upper obstacle $%
U$ is a quadruple of processes $(y,z,A,K)$ satisfying \newline
(\textbf{i}) $(y,z,A,K)\in \mathbf{D}_{\mathcal{F}}^{2}(0,T)\times \mathbf{L}%
_{\mathcal{F}}^{2}(0,T;\mathbb{R}^{d})\times (\mathbf{A}_{\mathcal{F}%
}^{2}(0,T))^{2}$ with $dA\geq 0$ and $dK\geq 0$ verifies
\begin{eqnarray}
y_{t}
&=&X+\int_{t}^{T}g(s,y_{s},z_{s})ds+A_{T}-A_{t}-(K_{T}-K_{t})-%
\int_{t}^{T}z_{s}dB_{s},  \label{RCBSDE} \\
&&\;
\begin{array}{c}
d_{\Gamma _{t}}(z_{t})=0,\;\;dP\times dt\mbox{-a.s. }\;\mathcal{V}_{[0,T]}[A-K]=%
\mathcal{V}_{[0,T]}[A+K],
\end{array}
\;  \nonumber
\end{eqnarray}
where $\mathcal{V}_{[0,T]}(\varphi )$ denotes the total variation of
a process $\varphi $ on $[0,T]$. \newline (\textbf{ii}) $y_{t}\leq
U_{t},$ $dP\times dt$-a.s., the generalized Skorohod reflecting
condition is satisfied:
\[
\int_{0}^{T}(U_{t-}^{*}-y_{t-})dK_{t}=0,\;\mbox{a.s., for any }
U^{*}\in \mathbf{D}_{\mathcal{F}}^{2}(0,T),\;\mbox{s.t.\ }y_{t}\geq
U_{t}^{*}\geq U_{t},\;dP\times dt\;\mbox{-a.s..}
\]
\textbf{(iii) }For any other quadruple $(y^{*},z^{*},A^{*},K^{*})$
satisfying \textbf{(i) }and \textbf{(ii)}, we have
\[
y_{t}\leq y_{t}^{*},\;\;\;\;\;0\leq t\leq T,\;\;\mbox{a.s. }
\]
\end{definition}

Like increasing processes of the solution of $g_{\Gamma}$-reflecting solution with one lower barrier,
here increasing processes $A$ and $K$ function separately. The role of $dA$ is to keep $z_{t}$ staying in the domain $\Gamma _{t}$, and $%
dK$ increases only when process $y_{t}$ tends to cross upwards the
upper barrier $U$.

\begin{remark}
The formula  $\mathcal{V}_{[0,T]}[A-K]=\mathcal{V}%
_{[0,T]}[A+K]$ in (\ref{RCBSDE}) implies that they never act at same
time. This helps us to separate two increasing processes completely.
And the proof of theorem \ref{existru} in subsection 3.3 shows that
$A$ and $K$ are just the limit of the corresponding terms in
penalization equations.
\end{remark}

Then we have the existence of the $g_{\Gamma }$--reflected solution with an upper obstacle $%
U$:

\begin{theorem}
\label{existru}Assume that (\ref{Lip}) holds for $g$ and (\ref{Gamma}) holds%
 for the constraint $\Gamma $, $U\ $is a $\mathcal{F}_{t}$-adapted
RCLL process satisfying (\ref{barc}). For a given terminal condition
$X\in
\mathbf{L}^{2}(\mathcal{F}_{T})$, the $g_{\Gamma }$--reflected solution $%
(y,z,A,K)$ with upper obstacle $U$  of Definition \ref{RCBSDE1}
\textbf{(i)-(iii)} exists.
\end{theorem}

\begin{remark}
For general case when $\Gamma $ depends on $y$, satisfying
(\ref{Gamma}), theorem \ref{existru} holds under the assumption of
the
existence of a special solution, i.e. there exists a quadruple $%
(y^{*},z^{*},A^{*},K^{*})\in \mathbf{D}_{\mathcal{F}}^{2}(0,T)\times \mathbf{%
L}_{\mathcal{F}}^{2}(0,T;\mathbb{R}^{d})\times (\mathbf{A}_{\mathcal{F}%
}^{2}(0,T))^{2}$, s.t. $dA_{t}^{*}\geq 0$, $dK_{t}^{*}\geq 0$ and
\begin{eqnarray}
y_{t}^{*}
&=&X+%
\int_{t}^{T}g(s,y_{s}^{*},z_{s}^{*})ds+(A_{T}^{*}-A_{t}^{*})-(K_{T}^{*}-K_{t}^{*})-\int_{t}^{T}z_{s}^{*}dB_{s},
\label{H1'} \\
&&\ \ \ \
\begin{array}{c}
d_{\Gamma _{t}}(y_{t}^{*},z_{t}^{*})=0\mbox{, a.s. a.e.}
\end{array}
\nonumber \\
y_{t}^{*} &\leq &U_{t}\mbox{, \ \ \ \ }%
\int_{0}^{T}(y_{t-}^{*}-U_{t-})dK_{t}^{*}=0\mbox{, a.s..}  \nonumber
\end{eqnarray}
This assumption is not easy to verify for general case. While if
$\Gamma _{t}=[L_{t},+\infty )$, it turns out to be a reflected BSDE
with two barriers $L$ and $U$, then refer to \cite{PX2005}, we know
that assumption (\ref{H1'}) can be changed to another sufficient
condition: there exists a semimartingale $X$, such that $L\leq X\leq
U$, $P$-a.s. a.e., which guarantee the existence of a special
solution.
\end{remark}

The proofs of Theorem \ref{existrl} and Theorem \ref{existru} are
given in the following subsections.

\subsection{Existence of $g_{\Gamma }$-reflected BSDE with a lower barrier:
Proof of Theorem \ref{existrl}}

We prove theorem \ref{existrl} by an approximation procedure. For
$m$, $n\in \mathbb{N}$, we consider the penalization equations,
\begin{eqnarray}
y_{t}^{m,n}
&=&X+\int_{t}^{T}g(s,y_{s}^{m,n},z_{s}^{m,n})ds+m\int_{t}^{T}d_{\Gamma
_{s}}(y_{s}^{m,n},z_{s}^{m,n})ds  \label{perefl} \\
&&\
+n\int_{t}^{T}(L_{s}-y_{s}^{m,n})^{+}ds-\int_{t}^{T}z_{s}^{m,n}dB_{s}.
\nonumber
\end{eqnarray}
Define $A_{t}^{m,n}=m\int_{0}^{t}d_{\Gamma
_{s}}(y_{s}^{m,n},z_{s}^{m,n})ds$ and
$\overline{A}_{t}^{m,n}=n\int_{0}^{t}(L_{s}-y_{s}^{m,n})^{+}ds$. We
have the following estimate.

\begin{lemma}
\label{est-rcbl}There exists a constant $C\in \mathbb{R}$ independent of $m$
and $n$, such that
\begin{equation}
E[\sup_{0\leq t\leq T}(y_{t}^{m,n})^{2}]+E\int_{0}^{T}\left|
z_{s}^{m,n}\right| ^{2}ds+E[(A_{T}^{m,n}+\overline{A}_{T}^{m,n})^{2}]\leq C.
\label{est-ryza}
\end{equation}
\end{lemma}

\noindent\textbf{Proof. } Set $m=n=0$, then we get a classical BSDE
\[
y_{t}^{0,0}=X+\int_{t}^{T}g(s,y_{s}^{0,0},z_{s}^{0,0})ds-%
\int_{t}^{T}z_{s}^{0,0}dB_{s}.
\]
For $(y^{*},z^{*},A^{*})$ given in (\ref{H1}), we have $d_{\Gamma
_{s}}(y_{s}^{*},z_{s}^{*})\equiv 0$ and $(L_{s}-y_{s}^{*})^{+}\equiv
0$, thus
\begin{eqnarray*}
y_{t}^{*} &=&X+\int_{t}^{T}g(s,y_{s}^{*},z_{s}^{*})ds+m\int_{t}^{T}d_{\Gamma
_{s}}(y_{s}^{*},z_{s}^{*})ds+n\int_{t}^{T}(L_{s}-y_{s}^{*})^{+}ds \\
&&+(A_{T}^{*}-A_{t}^{*})-%
\int_{t}^{T}z_{s}^{*}dB_{s},
\end{eqnarray*}
By comparison theorem, it follows $y_{t}^{*}\geq y_{t}^{m,n}\geq y_{t}^{0,0}$%
, $0\leq t\leq T$. So we have for some constant $C$ independent of $m$ and $%
n $,
\begin{equation}
E[\sup_{0\leq t\leq T}(y_{t}^{m,n})^{2}]\leq \max \{E[\sup_{0\leq t\leq
T}(y_{t}^{*})^{2}],\;E[\sup_{0\leq t\leq T}(y_{t}^{0,0})^{2}]\}\leq C.
\label{est-ry}
\end{equation}
Then applying It\^{o}'s formula to $\left| y_{t}^{m,n}\right| ^{2}$
and taking expectation, we get
\begin{eqnarray*}
&&\ \ \ E[\left| y_{t}^{m,n}\right| ^{2}]+E[\int_{t}^{T}\left|
z_{s}^{m,n}\right| ^{2}ds] \\
\ &\leq &E[X^{2}]+E\int_{t}^{T}g^{2}(s,0,0)ds+(2\mu +\mu
^{2})\int_{t}^{T}\left| y_{s}^{m,n}\right| ^{2}ds+\frac{1}{2}%
E[\int_{t}^{T}\left| z_{s}^{m,n}\right| ^{2}ds] \\
&&\ \ \ +\frac{1}{\alpha }E[\sup_{0\leq t\leq T}(y_{t}^{m,n})^{2}]+\alpha
E[(A_{T}^{m,n}-A_{t}^{m,n}+\overline{A}_{T}^{m,n}-\overline{A}%
_{t}^{m,n})^{2}],
\end{eqnarray*}
where $\alpha \in \mathbb{R}$ to be chosen later. Since $A_{t}^{m,n}$ and $%
\overline{A}_{t}^{m,n}$ are increasing processes, so
\begin{equation}
E\int_{0}^{T}\left| z_{s}^{m,n}\right| ^{2}ds\leq C+\alpha E[(A_{T}^{m,n}+%
\overline{A}_{T}^{m,n})^{2}].  \label{est-rz}
\end{equation}
While rewrite (\ref{perefl}) in the following form
\[
A_{T}^{m,n}+\overline{A}_{T}^{m,n}=y_{0}^{m,n}-X-%
\int_{0}^{T}g(s,y_{s}^{m,n},z_{s}^{m,n})ds+\int_{0}^{T}z_{s}^{m,n}dB_{s},
\]
then take square and expectation on both sides, we get
\begin{eqnarray*}
E[(A_{T}^{m,n}+\overline{A}_{T}^{m,n})^{2}] &\leq
&4E[(y_{0}^{m,n})^{2}]+4E[X^{2}]+16TE\int_{0}^{T}g^{2}(s,0,0)ds \\
&&\ \ \ +16\mu ^{2}TE\int_{0}^{T}\left| y_{s}^{m,n}\right| ^{2}ds+(16\mu
^{2}T+4)E\int_{0}^{T}\left| z_{s}^{m,n}\right| ^{2}ds,
\end{eqnarray*}
we then have
\begin{equation}
E[(A_{T}^{m,n}+\overline{A}_{T}^{m,n})^{2}]\leq C+(16\mu
^{2}T+4)E\int_{0}^{T}\left| z_{s}^{m,n}\right| ^{2}ds.  \label{est-ra}
\end{equation}
Compare (\ref{est-rz}) and (\ref{est-ra}), set $\alpha =\frac{1}{32\mu
^{2}T+8}$, we deduce (\ref{est-ryza}). \hfill $\Box $\medskip

\noindent \textbf{Proof of Theorem \ref{existrl}. }In
(\ref{perefl}), we fix $m\in \mathbb{N}$, and set
\[
g^{m}(t,y,z):=(g+md_{\Gamma _{t}})(t,y,z).
\]
This is a Lipschitz function. It follows from theorem 4.1 in \cite{PX2005}
that, as $n\rightarrow \infty $, with (\ref{est-ryza}) the triple $%
(y^{m,n},z^{m,n},\overline{A}^{m,n})$ converges to $(y^{m},z^{m},\overline{A}%
^{m})\in \mathbf{D}_{\mathcal{F}}^{2}(0,T)\times \mathbf{L}_{\mathcal{F}%
}^{2}(0,T)\times \mathbf{A}_{\mathcal{F}}^{2}(0,T)$, which is the
solution of the following reflected BSDE whose coefficient is
$g^{m}$:
\begin{eqnarray}
y_{t}^{m} &=&X+\int_{t}^{T}(g+md_{\Gamma _{s}})(s,y_{s}^{m},z_{s}^{m})ds+%
\overline{A}_{T}^{m}-\overline{A}_{t}^{m}-\int_{t}^{T}z_{s}^{m}dB_{s},
\label{pen-ref} \\
y_{t}^{m} &\geq &L_{t}\mbox{, \ \ \ }\int_{0}^{T}(y_{t-}-L_{t-}^{*})d%
\overline{A}_{t}^{m}=0,  \nonumber \\
&&
\begin{array}{c}
\mbox{for each }L^{*}\in \mathbf{D}_{\mathcal{F}}^{2}(0,T)\mbox{, such that }%
y\geq L^{*}\geq L,dP\times dt\mbox{ a.s.}.
\end{array}
\;  \nonumber
\end{eqnarray}
We denote $A_{t}^{m}=m\int_{0}^{t}d_{\Gamma _{s}}(z_{s}^{m})ds$. By (\ref
{est-ryza}) we have the following estimate:
\[
E[\sup_{0\leq t\leq T}(y_{t}^{m})^{2}]+E\int_{0}^{T}\left| z_{s}^{m}\right|
^{2}ds+E[(A_{T}^{m}+\overline{A}_{T}^{m})^{2}]\leq C.
\]
Then by comparison theorem \ref{CompR} for reflected BSDEs, we have $%
y_{t}^{m}\leq y_{t}^{m+1}$, $\overline{A}_{t}^{m}\geq \overline{A}_{t}^{m+1}$
and $d\overline{A}_{t}^{m}\geq d\overline{A}_{t}^{m+1}$ on $[0,T]$. Thus,
when $m\rightarrow \infty $, $y_{t}^{m}\nearrow y_{t}\leq y_{t}^{*}$,\ $%
\overline{A}_{t}^{m}\searrow \overline{A}_{t}$ in $\mathbf{L}^{2}(\mathcal{F}%
_{t})$, for each $t\in [0,T]$. Thanks to Fatou's lemma, we get $%
E[\sup_{0\leq t\leq T}\left| y_{t}\right| ^{2}]<\infty $, and thus $%
y^{m}\rightarrow y$ in $\mathbf{L}_{\mathcal{F}}^{2}(0,T)$ in view
of dominate convergence theorem. Since $\overline{A}^{m}$ is RCLL,
we can not directly apply the monotonic limit theorem, Theorem 2.1
in \cite{P99}. However it is easy to know that the limit $y$ can be
written in the following form
\[
y_{t}=y_{0}-\int_{0}^{t}g_{s}^{0}ds-A_{t}-\overline{A}_{t}+%
\int_{0}^{t}z_{s}dB_{s},
\]
where $z$ and $g^{0}$ (resp. $A_t$) are the weak limit of $z^{m}$
and $g^{m}$
(resp. $A^{m}_t$) in $\mathbf{L}_{\mathcal{F}}^{2}(0,T)$ (resp. $\mathbf{L}%
^{2}(\mathcal{F}_{t})$). By Lemma 2.2 in \cite{P99}, we know that
$y$ is RCLL. We then apply It\^{o}'s rule to $\left|
y_{t}^{m}-y_{t}\right| ^{2}$ on interval $[\sigma ,\tau ]$, with
stopping times $0\leq \sigma \leq \tau \leq T$. It follows that
\begin{eqnarray*}
&&\ \ E|y_{\sigma }^{m}-y_{\sigma }|^{2}+E\int_{\sigma
}^{\tau }|z_{s}^{m}-z_{s}|^{2}ds \\
\ &=&E|y_{\tau }^{m}-y_{\tau }|^{2}+E\sum_{t\in (\sigma
,\tau ]}[(\Delta A_{t})^{2}-(\overline{A}_{t}^{m}-\overline{A}_{t})^{2}]-2%
E\int_{\sigma }^{\tau }(y_{s}^{m}-y_{s})(g_{s}^{m}-g_{s}^{0})ds \\
&&\ \ +2E\int_{(\sigma ,\tau
]}(y_{s}^{m}-y_{s})dA_{s}^{m}-2E\int_{(\sigma ,\tau
]}(y_{s}^{m}-y_{s})dA_{s}+2E\int_{(\sigma ,\tau
]}(y_{s-}^{m}-y_{s-})d(\overline{A}_{s}^{m}-\overline{A}_{s}).
\end{eqnarray*}
Since $E\int_{(\sigma ,\tau ]}(y_{s}^{m}-y_{s})dA_{s}^{m}\leq 0$
and $E\int_{(\sigma ,\tau ]}(y_{s-}^{m}-y_{s-})d(\overline{A}%
_{s}^{m}-\overline{A}_{s})\leq 0$, so we get
\begin{eqnarray*}
E\int_{\sigma }^{\tau }|z_{s}^{m}-z_{s}|^{2}ds &\leq &E%
|y_{\tau }^{m}-y_{\tau }|^{2}+E\sum_{t\in (\sigma ,\tau ]}(\Delta
A_{t})^{2}+2E\int_{\sigma }^{\tau }\left| y_{s}^{m}-y_{s}\right|
\left| g_{s}^{m}-g_{s}^{0}\right| ds \\
&&\ +2E\int_{(\sigma ,\tau ]}\left| y_{s}^{m}-y_{s}\right| dA_{s}.
\end{eqnarray*}
Now we are in the same situation as in the proof of the monotonic limit
theorem (cf. \cite{P99}, Proof of Theorem 2.1). We then can follow the proof
and get $z^{m}\rightarrow z$ strongly in $\mathbf{L}_{\mathcal{F}}^{p}(0,T)$%
, for $p<2$.

From the Lipschitz property of $g$, we deduce that $(y,z,A,\overline{A})$
verify the equation
\[
y_{t}=X+\int_{t}^{T}g(s,y_{s},z_{s})ds+A_{T}-A_{t}+\overline{A}_{T}-%
\overline{A}_{t}-\int_{t}^{T}z_{s}dB_{s}.
\]
The estimate $E[(A_{T}^{m})^{2}]\leq C$ implies
$E[(\int_{0}^{T}d_{\Gamma _{s}}(z_{s}^{m})ds)^{2}]\leq
\frac{C}{m^{2}}$, thus
\[
E[\int_{0}^{T}d_{\Gamma _{s}}(z_{s})ds]=0,\;\;\mbox{or \ }d_{\Gamma
_{t}}(z_{t})\equiv 0\mbox{, }dP\times dt-a.s..
\]

It remains to prove that $(y,A)$ satisfies condition \textbf{(ii)} of
Definition \ref{RCBSDE2}, i.e., $y\geq L$ and
\begin{equation}
\mbox{ }\int_{0}^{T}(y_{t-}-L_{t-}^{*})d\overline{A}_{t}=0.
\label{GammaRefl}
\end{equation}
By $y^{m}\geq L$ we have $y\geq L$ and, for each $L^{*}\in \mathbf{D}_{%
\mathcal{F}}^{2}(0,T)\;$such that\ $y\geq L^{*}\geq L$,
\begin{eqnarray}
\int_{0}^{T}(y_{t-}-y_{t-}^{m}\wedge L_{t-}^{*})d\overline{A}_{t}
=&&\int_{0}^{T}(y_{t-}-y_{t-}^{m})d\overline{A}_{t}+%
\int_{0}^{T}(y_{t-}^{m}-y_{t-}^{m}\wedge L_{t-}^{*})d\overline{A}_{t}^{m}
\nonumber \\
& +&\int_{0}^{T}(y_{t-}^{m}-y_{t-}^{m}\wedge L_{t-}^{*})d(A_{t}-\overline{A}%
_{t}^{m}).  \nonumber
\end{eqnarray}
As $m\rightarrow \infty $, the first term on the right side tends to zero
due to Lebesgue domination theorem. The second term is null because of (\ref
{pen-ref}) and since $y^{m}\geq y^{m}\wedge L^{*}\geq L$. For the third term
we have
\begin{eqnarray}
E|\int_{0}^{T}(y_{t-}^{m}-y_{t-}^{m}\wedge L_{t-}^{*})d(A_{t}-\overline{A}%
_{t}^{m})|& \leq &E[\sup_{t\in [0,T]}|y_{t}^{m}-y_{t-}^{m}\wedge
L_{t-}^{*}|(A_{T}^{m}-A_{T})]  \nonumber \\
& \leq &E[\sup_{t\in [0,T]}|y_{t}^{m}-y_{t-}^{m}\wedge
L_{t-}^{*}|^{2}]^{1/2}E[(A_{T}^{m}-A_{T})^{2}]^{1/2}  \nonumber
\end{eqnarray}
which converges also to zero since $E[(A_{T}^{m}-A_{T})^{2}]^{1/2}\searrow 0$%
. Thus the left hand term must tend to zero. This with $y^{m}\wedge
L^{*}\nearrow L^{*}$ yields (\ref{GammaRefl}).

We now prove \textbf{(iii)}. Consider a quadruple
$(y^{*},z^{*},A^{*},\bar A^{*})$ which satisfies \textbf{(i)} and
\textbf{(ii)}. Since $d_{\Gamma _s}(y_s^{*},z_s^{*})\equiv 0$, we
have
\[
y_t^{*}=X+\int_t^Tg(s,y_s^{*},z_s^{*})ds+m\int_t^Td_{\Gamma
_s}(y_s^{*},z_s^{*})ds+A_T^{*}-A_t^{*}+\overline{A}_T^{*}-\overline{A}%
_t^{*}-\int_t^Tz_sdB_s.
\]
By comparison theorem \ref{CompR} it follows that $y^{*}\geq y^m$, for all $%
m $. Thus (iii) holds. \hfill$\Box$\medskip

\begin{remark}
\label{contil}If $L$ is continuous or only has positive jumps ($L_{t-}\leq
L_{t}$), then $\overline{A}$ is a continuous process. In this case, in (\ref
{penrbsde}), $\overline{A}^{n}$ are continuous, and $\overline{A}%
_{t}^{n}\geq \overline{A}_{t}^{n+1}$, $d\overline{A}_{t}^{n}\geq d\overline{A%
}_{t}^{n+1}$, $0\leq t\leq T$, with $E[(\overline{A}_{T}^{n})^{2}]\leq C$.
Then $\overline{A}_{t}^{n}\searrow \overline{A}_{t}$, $0\leq t\leq T$.
Moreover
\[
0\leq \overline{A}_{t}^{n}-\overline{A}_{t}\leq \overline{A}_{T}^{n}-%
\overline{A}_{T}.
\]
Thus we have uniform convergence:
\[
E[\sup_{0\leq t\leq T}(\overline{A}_{t}^{n}-\overline{A}_{t})^{2}]\leq E[(%
\overline{A}_{T}^{n}-\overline{A}_{T})^{2}]\rightarrow 0,\;\;\mbox{as \ }%
n\rightarrow \infty .
\]
\end{remark}

\subsection{Some convergence results of $g_{\Gamma }$-reflected solution with a lower barrier}

As we know, the reflected BSDE can be considered as a special kind
of constraint BSDE, with $\Gamma _{t}=[L_{t},+\infty )\times
\mathbb{R}$. If we put two constraint together, i.e. set
$\widehat{\Gamma }_{t}=\Gamma _{t}\cap [L_{t},+\infty )$, then the
penalization equation becomes the following one: for $n\in
\mathbb{N}$
\begin{eqnarray}
y_{t}^{n,n} &=&X+\int_{t}^{T}g(s,y_{s}^{n,n},z_{s}^{n,n})ds+n\int_{t}^{T}d_{%
\widehat{\Gamma }_{s}}(y_{s}^{n,n},z_{s}^{n,n})ds-%
\int_{t}^{T}z_{s}^{n,n}dB_{s}  \label{Pen-nn} \\
&=&X+\int_{t}^{T}g(s,y_{s}^{n,n},z_{s}^{n,n})ds+n\int_{t}^{T}d_{\Gamma
_{s}}(y_{s}^{n,n},z_{s}^{n,n})ds+n\int_{t}^{T}(L_{s}-y_{s}^{n,n})^{+}ds
\nonumber \\
&&-\int_{t}^{T}z_{s}^{n,n}dB_{s}.  \nonumber
\end{eqnarray}
Setting $\widehat{A}_{t}^{n}=n\int_{0}^{t}d_{\widehat{\Gamma }%
_{s}}(y_{s}^{n,n},z_{s}^{n,n})ds$, with monotonic limit theorem in
\cite{P99}, we know that let $n\rightarrow \infty $,
$(y^{n,n},z^{n,n},\widehat{A}^{n,n})$
converges to $(\widehat{y},\widehat{z},\widehat{A})\in \mathbf{L}_{\mathcal{F%
}}^{2}(0,T)\times \mathbf{L}_{\mathcal{F}}^{2}(0,T;\mathbb{R}^{d})\times
\mathbf{A}_{\mathcal{F}}^{2}(0,T)$, where
\[
\widehat{y}_{t}=X+\int_{t}^{T}g(s,\widehat{y}_{s},\widehat{z}_{s})ds+%
\widehat{A}_{T}-\widehat{A}_{t}-\int_{t}^{T}\widehat{z}_{s}dB_{s}.
\]
Then we have

\begin{proposition}
\label{conv-nn}The two limits are equal in the following sense:
\[y_{t}=%
\widehat{y}_{t}, z_{t}=\widehat{z}_{t}, A_{t}+\overline{A}_{t}=\widehat{A%
}_{t}.
\]
\end{proposition}

Before we give the proof of this proposition, we consider another
way to prove the convergence by the penalization equations given by
(\ref{perefl}), i.e. first let $m\rightarrow \infty $, then let
$n\rightarrow \infty $,
while in former subsection, we get the $g_{\Gamma }$-reflected solution $%
(y,z,A,\overline{A})$ of Definition \ref{dif-RCBSDE2}, by first letting $%
n\rightarrow \infty $, then letting $m\rightarrow \infty $. So as $%
m\rightarrow \infty $, we get that the triple $(y^{m,n},z^{m,n},A^{m,n})$
converges to $(y^{n},z^{n},A^{n})\in \mathbf{D}_{\mathcal{F}}^{2}(0,T)\times
\mathbf{L}_{\mathcal{F}}^{2}(0,T;\mathbb{R}^{d})\times \mathbf{A}_{\mathcal{F%
}}^{2}(0,T)$, which is the solution of constraint BSDE with coefficient $%
g^{n}=g+n(L_{t}-y)^{+}$:
\begin{eqnarray}
y_{t}^{n}
&=&X+\int_{t}^{T}g(s,y_{s}^{n},z_{s}^{n})ds+A_{T}^{n}-A_{t}^{n}+n%
\int_{t}^{T}(L_{s}-y_{s}^{n})^{+}ds-\int_{t}^{T}z_{s}^{n}dB_{s},
\label{pen-n} \\
(z_{t}^{n}) &\in &\Gamma _{t},\;\;dP\times dt \mbox{-a.s.,
}\;\;dA^{n}\geq 0. \nonumber
\end{eqnarray}
Define $\overline{A}_{t}^{n}=n\int_{0}^{t}(L_{s}-y_{s}^{n})^{+}ds$. With
same method in former subsection, we can prove that as $n\rightarrow \infty $%
, $(y^{n},z^{n},A^{n},\overline{A}^{n})$ converges to $(\widetilde{y},%
\widetilde{z},\widetilde{A},\widetilde{\overline{A}},)$ where
\[
\widetilde{y}_{t}=X+\int_{t}^{T}g(s,\widetilde{y}_{s},\widetilde{z}_{s})ds+%
\widetilde{A}_{T}-\widetilde{A}_{t}+\widetilde{\overline{A}}_{T}-\widetilde{%
\overline{A}}_{t}-\int_{t}^{T}\widetilde{z}_{s}dB_{s}.
\]
Then we have

\begin{proposition}
\label{2limits}The two limits are equal, in the following sense,
\[y_{t}=\widetilde{y}%
_{t}, z_{t}=\widetilde{z}_{t} \mbox{ and } A_{t}+\overline{A}_{t}=\widetilde{A}%
_{t}+\widetilde{\overline{A}}_{t}, 0\leq t\leq T.
\]
\end{proposition}

\noindent\textbf{Proof. }
By comparison theorem for (\ref{perefl}) and (\ref{pen-ref}), we have $%
y_{t}^{m,n}\leq y_{t}^{m}$, which follows $y_{t}^{n}\leq y_{t}$, when
letting $m\rightarrow \infty $. Then let $n\rightarrow \infty $, we get $\widetilde{y}%
_{t}\leq y_{t}$. Symmetrically compare (\ref{perefl}) and (\ref{pen-n}), $%
y_{t}^{m,n}\leq y_{t}^{n}$, let $n\rightarrow \infty $, we get $%
y_{t}^{m}\leq \widetilde{y}_{t}$, then as $m\rightarrow \infty $, it
follows $y_{t}\leq \widetilde{y}_{t}$. So $y_{t}=\widetilde{y}_{t}$,
$0\leq t\leq T$. The rest follows easily. $\square $

Now we prove proposition \ref{conv-nn}:

\textbf{Proof of proposition \ref{conv-nn}}: For $m\leq n$, by comparison
theorem for (\ref{perefl}) and (\ref{Pen-nn}), we have $y_{t}^{m,n}\leq
y_{t}^{n,n}$. Let $n\rightarrow \infty $, then $m\rightarrow \infty $, we
get
\[
y_{t}\leq \widehat{y}_{t}.
\]
Similarly, for $m\geq n$, using again comparison theorem, we have $%
y_{t}^{m,n}\geq y_{t}^{n,n}$. First let $m\rightarrow \infty $, then $%
n\rightarrow \infty $, it follows
\[
\widetilde{y}_{t}\geq \widehat{y}_{t}.
\]
With proposition \ref{2limits}, we obtain $y_{t}=\widetilde{y}_{t}=\widehat{y%
}_{t}$. Other equalities follow easily. $\square $

These results show that for $g_{\Gamma }$-reflected BSDE with a
lower barrier, we can get its solution via penalisation equations by
different convergence method. No matter letting $m\rightarrow \infty
$ first or
letting $n\rightarrow \infty $ first, even considering dialogue sequence $%
(m=n)$, the limits we get are the same. By (\ref{Pen-nn}) and
monotonic limit theorem in \cite{P99}, we get $g_{\widehat{\Gamma
}}$-solution $\widehat{y}$
directly, increasing process $\widehat{A}$ is to keep $(y,z)$ stay in $%
\widehat{\Gamma} $, but we do not know any further property. But the $g_{\Gamma }$%
-reflected solution, i.e. definition \ref{RCBSDE2}, permits us to
have a
decomposition of $\widehat{A}$, with $\widehat{A}=A+\overline{A}$, where $%
\overline{A}$ serves for $y_{t}$ to get $y_{t}\geq L_{t}$ and $A$
serves for $z_{t}$ to keep $z_{t}\in \Gamma _{t}$, $dP\times
dt$-a.s.. And this property plays an important role when we study
the American option in incomplete market.

\begin{remark}
Proposition \ref{conv-nn} is still true if we consider the more
general case $\Gamma $ could depend on $y$, which satisfies
(\ref{Gamma}). Moreover we can generalize the constraint of
reflecting with a lower barrier $L$ by another general constraint
$\Lambda (t,\omega )$ which satisfies (\ref{Gamma}), and Proposition
\ref{conv-nn} still holds.
\end{remark}

\subsection{Existence of $g_{\Gamma }$-reflected solution with an upper barrier:
Proof of Theorem \ref{existru}}

For each $n\in \mathbb{N}$, we consider the solution
$
(y^{n},z^{n},K^{n})\in \mathbf{D}_{\mathcal{F}}^{2}(0,T)\times \mathbf{L}_{%
\mathcal{F}}^{2}(0,T;\mathbb{R}^{d})\times \mathbf{A}_{\mathcal{F}}^{2}(0,T)
$
of the following reflected BSDE with the coefficient $g^{n}=(g+nd_{\Gamma
_{t}})(t,y,z)$ and the upper reflecting obstacle $U$:
\begin{eqnarray}
y_{t}^{n} &=&X+\int_{t}^{T}(g+nd_{\Gamma
_{t}})(s,y_{s}^{n},z_{s}^{n})ds-(K_{T}^{n}-K_{t}^{n})-%
\int_{t}^{T}z_{s}^{n}dB_{s},  \label{penrbsde} \\
y^{n} &\leq &U,\;dP\times dt\mbox{-a.s.}\;\;dK\geq 0\mbox{,\ \ \ and \ }%
\int_{0}^{T}(U_{t-}^{*}-y_{t-}^{n})dK_{t}^{n}=0,  \nonumber \\
\forall \;U^{*} &\in &\mathbf{D}_{\mathcal{F}}^{2}(0,T),\;\;\mbox{such that \ }%
y^{n}\leq U^{*}\leq U\;\;dP\times dt\mbox{-a.s.}.  \nonumber
\end{eqnarray}
Since $g^{n}$ is Lipschitz with respect to $(y,z)$, this equation has a
unique solution. We denote
\[
A_{t}^{n}=n\int_{0}^{t}d_{\Gamma _{s}}(y_{s}^{n},z_{s}^{n})ds.
\]

Before to prove the a priori estimation for $(y^{n},z^{n},A^{n},K^{n})$, we
need the following lemma.

\begin{lemma}
\label{special}For $X\in \mathbf{L}^{2}(\mathcal{F}_{T})$, there exists a
quadruple $(y^{*},z^{*},A^{*},K^{*})\in \mathbf{D}_{\mathcal{F}%
}^{2}(0,T)\times \mathbf{L}_{\mathcal{F}}^{2}(0,T;\mathbb{R}^{d})\times (%
\mathbf{A}_{\mathcal{F}}^{2}(0,T))^{2}$ satisfies
\begin{eqnarray}
y_{t}^{*}
&=&X+%
\int_{t}^{T}g(s,y_{s}^{*},z_{s}^{*})ds+(A_{T}^{*}-A_{t}^{*})-(K_{T}^{*}-K_{t}^{*})-\int_{t}^{T}z_{s}^{*}dB_{s},
\label{cond-refl} \\
\ d_{\Gamma _{t}}(z_{t}^{*}) &=&0\mbox{, }dP\times dt\mbox{-a.s. and
}\
y_{t}^{*}\leq U_{t}\mbox{, \ \ \ \ }%
\int_{0}^{T}(y_{t-}^{*}-U^*_{t-})dK_{t}^{*}=0\mbox{, a.s..}
\nonumber\\
\forall \;U^{*} &\in &\mathbf{D}_{\mathcal{F}}^{2}(0,T),\;\;\mbox{such that \ }%
y^{*}\leq U^{*}\leq U\;\;dP\times dt\mbox{-a.s.}.  \nonumber
\end{eqnarray}
\end{lemma}

\noindent\textbf{Proof. }
Fix a process $\sigma _{t}\in \mathbf{L}_{\mathcal{F}}^{2}(0,T;\mathbb{R}%
^{d})$ satisfying $\sigma _{t}\in \Gamma _{t}$, $t\in [0,T]$. We consider a
forward SDE with an upper barrier $U_{t}$%
\begin{eqnarray*}
dx_{t} &=&-g(t,x_{t},\sigma _{t})dt-d\overline{A}_{t}+\sigma _{t}dB_{t}, \\
x_{0} &=&1\wedge U_{0}.
\end{eqnarray*}
Here $\overline{A}$ is a process in
$\mathbf{A}_{\mathcal{F}}^{2}(0,T)$, such that $x_{t}\leq U_{t}$,
a.s. a.e.. Set
\[
y_{t}^{*}=x_{t},z_{t}^{*}=\sigma _{t},A_{t}^{*}=\overline{A}%
_{t}+(x_{T}-X)^{+}1_{\{t=T\}},K_{t}^{*}=(x_{T}-X)^{+}1_{\{t=T\}}.
\]
Then this quadruple is just the one we need. $\square $

\begin{lemma}
We have the following estimates: there exists a constant $C>0$, independent
of $n$, such that
\begin{equation}
E[\sup_{0\leq t\leq T}(y_{t}^{n})^{2}]+E\int_{0}^{T}\left| z_{s}^{n}\right|
^{2}ds+E[(A_{T}^{n})^{2}]+E[(K_{T}^{n})^{2}]\leq C.  \label{est-YK}
\end{equation}
\end{lemma}

\proof%
Consider the following reflected BSDE with $U$ as its upper reflecting
obstacle,
\begin{eqnarray*}
y_{t}^{0}
&=&Y_{T}+\int_{t}^{T}g(s,y_{s}^{0},z_{s}^{0})ds-(K_{T}^{0}-K_{t}^{0})-%
\int_{t}^{T}z_{s}^{0}dB_{s},\;\;t\in [0,T], \\
y_{t}^{0} &\leq &U_{t},\;dK_{t}\geq
0,\;\;\int_{0}^{T}(y_{t-}^{0}-U^*_{t-})dK_{t}^{0}=0.\\
\forall \;U^{*} &\in &\mathbf{D}_{\mathcal{F}}^{2}(0,T),\;\;\mbox{such that \ }%
y^{0}\leq U^{*}\leq U\;\;dP\times dt\mbox{-a.s.}.  \nonumber
\end{eqnarray*}
This equation has a unique solution $(y^{0},z^{0},K^{0})\in \mathbf{D}_{%
\mathcal{F}}^{2}(0,T)\times \mathbf{L}_{\mathcal{F}}^{2}(0,T;\mathbb{R}%
^{d})\times \mathbf{A}_{\mathcal{F}}^{2}(0,T)$. By comparison theorem of
reflected BSDEs $y^{n}\geq y^{0}$.

On the other hand, from proposition \ref{special}, there exists $%
(y^{*},z^{*},A^{*},K^{*}) $ satisfying
\begin{eqnarray*}
y_{t}^{*} &=&Y_{T}+\int_{t}^{T}(g+nd_{\Gamma
_{s}})(s,y_{s}^{*},z_{s}^{*})ds+(A_{T}^{*}-A_{t}^{*})-(K_{T}^{*}-K_{t}^{*})-%
\int_{t}^{T}z_{s}^{*}dB_{s}, \\
y_{t}^{*} &\leq &U_{t}\mbox{, \ \ }%
\int_{0}^{T}(y_{t-}^{*}-U_{t-})dK_{t}^{*}=0\mbox{, a.s..}
\end{eqnarray*}
It follows from the comparison theorem\ref{CompR} for reflected
BSDEs that
for each $n\in \mathbb{N}$, we have $y_{t}^{n}\leq y_{t}^{*}$, $%
K_{t}^{n}\leq K_{t}^{*}$ and $dK_{t}^{n}\leq dK_{t}^{*}$, $t\in [0,T].$ Thus
there exists a constant $C>0$, independent of $n$, such that
\begin{equation}
E[\sup_{0\leq t\leq T}(y_{t}^{n})^{2}]\leq E[\sup_{0\leq t\leq
T}\{(y_{t}^{0})^{2}+(y_{t}^{*})^{2}\}]\leq C.  \label{est-ys}
\end{equation}
and
\begin{equation}
E[(K_{T}^{n})^{2}]\leq E[(K_{T}^{*})^{2}]\leq C.  \label{est-k}
\end{equation}
To estimate $(z^{n},A^{n})$, we apply It\^{o}'s formula to $\left|
y_{t}^{n}\right| ^{2}$ then get
\begin{eqnarray*}
E[\left| y_{0}^{n}\right| ^{2}]+E[\int_{t}^{T}\left|
z_{s}^{n}\right| ^{2}ds] &\leq
&E[Y_{T}^{2}]+E\int_{0}^{T}(g(s,0,0))^{2}ds+(2\mu +\mu
^{2})\int_{0}^{T}\left| y_{s}^{n}\right| ^{2}ds+\frac{1}{2}%
E[\int_{0}^{T}\left| z_{s}^{n}\right| ^{2}ds] \\
&&+(\frac{1}{\alpha }+1)E[\sup_{0\leq t\leq
T}(y_{t}^{n})^{2}]+\alpha E[(A_{T}^{n})^{2}]+E[(K_{T}^{n})^{2}],
\end{eqnarray*}
where $\alpha $ is a positive constants to be chosen later. This
with the above two estimates (\ref{est-ys}) and (\ref{est-k}) yields
\begin{equation}
E[\int_{0}^{T}\left| z_{s}^{n}\right| ^{2}ds]\leq C+\alpha
E[(A_{T}^{n})^{2}].  \label{est-zas}
\end{equation}
On the other hand, again by (\ref{penrbsde}),
\begin{equation}
A_{T}^{n}=y_{0}^{n}-y_{T}^{n}-%
\int_{0}^{T}g(s,y_{s}^{n},z_{s}^{n})ds+K_{T}^{n}-\int_{0}^{T}z_{s}^{n}dB_{s}.
\label{est-zs}
\end{equation}
Thus
\begin{eqnarray*}
E[(A_{T}^{n})^{2}] &\leq
&5E[(y_{0}^{n})^{2}+(y_{T}^{n})^{2}+(K_{T}^{n})^{2}]+15TE%
\int_{0}^{T}(g(s,0,0))^{2}ds \\
&&\ \ \ \ \ +15\mu ^{2}TE\int_{0}^{T}\left| y_{s}^{n}\right| ^{2}ds+(15\mu
^{2}T+5)E\int_{0}^{T}\left| z_{s}^{n}\right| ^{2}ds,
\end{eqnarray*}
then
\begin{equation}
E[(A_{T}^{n})^{2}]\leq C+(15\mu ^{2}T+5)E\int_{0}^{T}\left| z_{s}^{n}\right|
^{2}ds.  \label{est-as}
\end{equation}
With (\ref{est-zas}), setting $%
\alpha =\frac{1}{30\mu ^{2}T+10}$, we finally obtain (\ref{est-YK}). \hfill $%
\Box $\medskip

\noindent \textbf{Proof of Theorem \ref{existru}.} In
(\ref{penrbsde}), since $g^{n}(t,y,z)\leq g^{n+1}(t,y,z)$, by
comparison theorem \ref{CompR} for reflected
BSDEs, $y^{0}\leq y^{n}\leq y^{n+1}\leq y^{*}$. Thus $%
\{y^{n}\}_{n=1}^{\infty }$ increasingly converges to $y$ as $n\rightarrow
\infty $, and
\[
E[\sup_{0\leq t\leq T}(y_{t})^{2}]\leq C.
\]
We also have
\[
\lim_{n\rightarrow \infty }E[\int_{0}^{T}\left| y_{t}^{n}-y_{t}\right|
^{2}dt]=0.
\]
Moreover from comparison theorem \ref{CompR}, we have $%
K_{t}^{n}\leq K_{t}^{n+1}\leq K_{t}^{*}$ and $dK_{t}^{n}\leq
dK_{t}^{n+1}\leq dK_{t}^{*}$, $0\leq t\leq T$. It follows that $%
\{K^{n}\}_{n=1}^{n}$ increasingly converges to an increasing process $K\in
\mathbf{A}_{\mathcal{F}}^{2}(0,T)$ with $E[(K_{T})^{2}]\leq C$. Moreover $%
A^{n}$ are continuous increasing processes with $E[(A_{T}^{n})^{2}]\leq C$.
From (\ref{est-zas}), there exists a process $z\in \mathbf{L}_{\mathcal{F}%
}^{2}(0,T;\mathbb{R}^{d})$, such that $z^{n}\rightarrow z$ weakly in $%
\mathbf{L}_{\mathcal{F}}^{2}(0,T;\mathbb{R}^{d})$.

Now the conditions of the generalized monotonic limit theorem,
Theorem 3.1 in \cite{PX2005}
are satisfied. Then we have $z^{n}\rightarrow z$ strongly in $\mathbf{L}_{%
\mathcal{F}}^{p}(0,T;\mathbb{R}^{d})$, for $p<2$. With the Lipschitz
condition of $g$, the limit $y\in \mathbf{D}^{2}(0,T)$ can be written as
\[
y_{t}=X+\int_{t}^{T}g(s,y_{s},z_{s})ds+(A_{T}-A_{t})-(K_{T}-K_{t})-%
\int_{t}^{T}z_{s}dB_{s},
\]
where, for each $t,$ $A_{t}^{n}\rightarrow A_{t}$ weakly in $\mathbf{L}^{2}(%
\mathcal{F}_{t})$, $K_{t}^{n}\rightarrow K_{t}$ strongly in $\mathbf{L}^{2}(%
\mathcal{F}_{t})$. $A,K\in \mathbf{A}_{\mathcal{F}}^{2}(0,T)$ are increasing
processes.

From $E[(A_T^n)^2]=E[(n\int_0^Td_{\Gamma _s}(y_s^n,z_s^n)ds)^2]\leq C$, it
follows that
\[
E[(\int_0^Td_{\Gamma _s}(y_s^n,z_s^n)ds)^2]\leq \frac C{n^2},
\]
while $d_{\Gamma _s}(y_s^n,z_s^n)\geq 0$, we get that $\int_0^Td_{\Gamma
_s}(y_s^n,z_s^n)ds\rightarrow 0$, as $n\rightarrow \infty $. With the
Lipschitz property of $d_{\Gamma _t}(y,z)$ and the convergence of $y^n$ and $%
z^n$, we deduce that
\[
d_{\Gamma _t}(y_t,z_t)=0\;\;dP\times dt\mbox{-a.s..}
\]

Now we consider \textbf{(ii)}. From $y^{n}\geq U$ we have $y\geq U$,
with
\[
\int_{0}^{T}(y_{t-}^{n}-U_{t-}^{*})dK_{t}^{n}=0,\;\;\forall U^{*}\in \mathbf{%
D}_{\mathcal{F}}^{2}(0,T)\mbox{, s.t. }U\geq U^{*}\geq y^{n}.
\]
Now, for each $U^{*}\in \mathbf{D}_{\mathcal{F}}^{2}(0,T)$, s.t. $U\geq
U^{*}\geq y$, since $y\geq y^{n}$, thus $U\geq U^{*}\geq y^{n}$
\[
\int_{0}^{T}(y_{t-}^{n}-U_{t-}^{*})dK_{t}^{n}=0\;\;\Rightarrow
\;\int_{0}^{T}(y_{t-}-U_{t-}^{*})dK_{t}^{n}=0.\;
\]
Recall that $dK_{t}^{n}\leq dK_{t}$, and $K_{T}^{n}\nearrow K_{T}$ in $%
\mathbf{L}^{2}(\mathcal{F}_{T})$, then
\[
0\leq \int_{0}^{T}(U_{t-}^{*}-y_{t-})d(K_{t}-K_{t}^{n})\leq \sup_{t\in
[0,T]}(U_{t}^{*}-y_{t})\cdot [K_{T}-K_{T}^{n}],
\]
and with the estimate of $y$ and (\ref{barc}), it follows
\textbf{(ii)} of Definition \ref{RCBSDE1} holds.

We now prove \textbf{(iii)}. In fact, for any other quadruple $(\overline{y},%
\overline{z},\overline{A},\overline{K})\in \mathbf{D}_{\mathcal{F}%
}^{2}(0,T)\times \mathbf{L}_{\mathcal{F}}^{2}(0,T;\mathbb{R}^{d})\times (%
\mathbf{A}^{2}(0,T))^{2}$ satisfying
\begin{eqnarray*}
\overline{y}_{t} &=&X+\int_{t}^{T}g(s,\overline{y}_{s},\overline{z}_{s})ds+%
\overline{A}_{T}-\overline{A}_{t}-(\overline{K}_{T}-\overline{K}%
_{t})-\int_{t}^{T}\overline{z}_{s}dB_{s}, \\
&&\
\begin{array}{c}
d_{\Gamma _{t}}(\overline{y}_{t},\overline{z}_{t})=0\mbox{, \ \ \ }d\bar{A}%
\geq 0,\;\;d\bar{K}\geq 0,
\end{array}
\\
\overline{y}_{t} &\leq &U_{t}^{*},\;\;\;\;\int_{0}^{T}(U_{t-}^{*}-\overline{y%
}_{t-})d\overline{K}_{t}=0\mbox{, \ \ a.s.,}
\end{eqnarray*}
for any $U^{*}\in \mathbf{D}_{\mathcal{F}}^{2}(0,T)$, such that \ $\overline{%
y}\leq U^{*}\leq U$. Then it also satisfies
\[
\overline{y}_{t}=X+\int_{t}^{T}g(s,\overline{y}_{s},\overline{z}%
_{s})ds+n\int_{t}^{T}d_{\Gamma _{s}}(\overline{y}_{s},\overline{z}_{s})ds+%
\overline{A}_{T}-\overline{A}_{t}-(\overline{K}_{T}-\overline{K}%
_{t})-\int_{t}^{T}\overline{z}_{s}dB_{s}.
\]
Compare it to (\ref{penrbsde}), we have $\overline{y}\geq y^{n}$, and $%
\overline{K}\geq K^{n}$. Let $n\rightarrow \infty $, it follows
\begin{equation}
\overline{y}_{t}\geq y_{t},\;\;\overline{K}\geq K_{t},\;\forall t\in [0,T],\;%
\mbox{a.s..}\;  \label{comsma}
\end{equation}
So $y$ is the smallest process satisfying Definition \ref{RCBSDE1}
\textbf{(i)} and \textbf{(ii)}.

It remains to prove the relation of the total variation in
(\ref{RCBSDE}) holds.
In fact, if it is not the case, set $\widetilde{V}_{t}=\mathcal{V}%
_{[0,t]}(A+K)$, then we define Jordan decomposition:
\[
\widetilde{A}_{t}=\frac{1}{2}(\widetilde{V}_{t}+A_{t}-K_{t}),\;\;\widetilde{K%
}_{t}=\frac{1}{2}(\widetilde{V}_{t}-A_{t}+K_{t}).
\]
With $d\widetilde{K}_{t}=\frac{1}{2}(d\widetilde{V}_{t}-dA_{t}+dK_{t})\leq
dK_{t}$ we have, for each $U^{*}\in \mathbf{D}_{\mathcal{F}}^{2}(0,T)$ with $%
U\geq U^{*}\geq y$,
\[
0\leq \int_{0}^{T}(y_{t-}-U_{t-}^{*})d\widetilde{K}_{t}\leq
\int_{0}^{T}(y_{t-}-U_{t-}^{*})dK_{t}=0.\;
\]
But in considering the second inequality of (\ref{comsma}), we have $%
\widetilde{K}\geq K$, which draws a contradiction. This completes
the proof. \hfill $\Box $\medskip

\begin{remark}
From the smallest property of $y-K$, it is the $g^{K}_{\Gamma
}$-solution with terminal condition $X-K_{T}$ , where
\[
g^{K}(t,y,z)=g(t,y+K_{t},z).
\]
\end{remark}

\begin{remark}
\label{contu}If $U$ is continuous (or satisfies $U_{t-}\geq U_{t}$),
then $K$ is a continuous process. In fact, by \cite{EKPPQ}, the
solution $y^{n}$ of (\ref{penrbsde}) as well as the reflecting
process $K^{n}$ are continuous. This with $K^{n}\leq K^{n+1}$ and
$dK^{n}\leq dK$ yields
\[
0\leq K_{t}-K_{t}^{n}\leq K_{T}-K_{T}^{n},
\]
and thus
\[
E[\sup_{0\leq t\leq T}(K_{t}-K_{t}^{n})^{2}]\leq
E[(K_{T}-K_{T}^{n})^{2}]\rightarrow 0.
\]
It follows that $K^{n}$ converges uniformly to $K$ on $[0,T]$. Thus
$K$ is continuous.
\end{remark}

\section{$g_{\Gamma }$-super(sub)martingales and its Doob-Meyer's type decomposition theorems}

Now we introduce the definitions of $g_{\Gamma }$-martingale, $g_{\Gamma }$%
-supermartingale and $g_{\Gamma }$-submartingale, by $g_{\Gamma }$%
-expectation introduced in section 3. Suppose that $g $ satisfies
(\ref{g-Gamma}), $g(t,0,0)=0$ and $(0,0)\in \Gamma _{t}$.

\begin{definition}
\label{smart}A process $Y\in \mathbf{D}_{\mathcal{F}}^{2}(0,T)$ is called a $%
g_{\Gamma }$-supermartingale (resp. $g_{\Gamma }$-submartingale) on
$[0,T]$, if for stopping times $\sigma $, $\tau $ valued in $[0,T]$,
with $\sigma
\leq \tau $, we have $Y_{\tau }\in \mathbf{L}_{+,\infty }^{2}(\mathcal{F}%
_{\tau })$ and
\[
\mathcal{E}_{\sigma ,\tau }^{g_{\Gamma }}[Y_{\tau }]\leq Y_{\sigma
},\mbox{ (resp. }\geq Y_{\sigma }\mbox{).}
\]

It is called a $g_{\Gamma }$--martingale if it is both a $g_{\Gamma }$%
--supermartingale and $g_{\Gamma }$--submartingale.
\end{definition}

The nonlineare Doob-Meyer's type decomposition theorem for $g$%
-super(sub)martingale in \cite{P99} plays un important role in theory of $g$%
-expectation. For $g_{\Gamma }$-super(sub)martingale, we have also
Doob-Meyer's type decomposition theorem. In fact, in \cite{PX2006},
we have proved the decomposition for $g_{\Gamma }$-supermartingale,
partly for $g_{\Gamma }$-submartingale. For completeness of this
paper, we still present the proofs. And these proofs are important
applications of $g_{\Gamma}$ reflecting solutions.

\subsection{$g_{\Gamma }$-supermartingale decomposition theorem}

In this section, we study the Doob-Meyer's type decomposition theorem for $%
g_{\Gamma }$-supermartingale. Before present the main result, we
first give a useful property of $g_{\Gamma }$-supermartingale:

\begin{proposition}
\label{gmphi}A process $Y$ is a $g_{\Gamma }$-supermartingale on $[0,T]$, if
and only if for all $m\geq 0$, it is a $(g+md_{\Gamma })$-supermartingale on
$[0,T]$.
\end{proposition}

\noindent\textbf{Proof. } We fix $t\in [0,T]$, and set
$y_{s}^{t}=\mathcal{E}_{s,t}^{g_{\Gamma
}}[Y_{t}]$, $0\leq s\leq t$. Let $(y^{t},z^{t},A^{t})$ be the $g_{\Gamma }$%
--solution on $[0,t]$:
\begin{eqnarray*}
y_{s}^{t}
&=&Y_{t}+\int_{s}^{t}g(r,y_{r}^{t},z_{r}^{t})dr+A_{t}^{t}-A_{s}^{t}-%
\int_{s}^{t}z_{r}^{t}dB_{r}, \\
&&
\begin{array}{c}
d_{\Gamma _{s}}(y_{s}^{t},z_{s}^{t})=0\mbox{, \ }s\in [0,t]\mbox{.}
\end{array}
\end{eqnarray*}
Consider the following penalization equation
\[
y_{s}^{t,m}=Y_{t}+\int_{s}^{t}g(r,y_{r}^{t,m},z_{r}^{t,m})dr+m%
\int_{s}^{t}d_{\Gamma
_{r}}(y_{r}^{t,m},z_{r}^{t,m})dr-\int_{s}^{t}z_{r}^{t,m}dB_{r},
\]
We observe that the above $(y^{t},z^{t},A^{t})$ also satisfies
\[
y_{s}^{t}=Y_{t}+\int_{s}^{t}g(r,y_{r}^{t},z_{r}^{t})dr+m\int_{s}^{t}d_{%
\Gamma
_{r}}(y_{r}^{t},z_{r}^{t})dr+A_{t}^{t}-A_{s}^{t}-%
\int_{s}^{t}z_{r}^{t}dB_{r}.
\]
From comparison theorem, we get $y^{t,m}\leq y^{t}$ on $[0,t]$. Thus
\[
\mathcal{E}_{s,t}^{g+md_{\Gamma }}[Y_{t}]\leq \mathcal{E}_{s,t}^{g_{\Gamma
}}[Y_{t}]\leq Y_{s},\;\;\forall m\geq 0.
\]
It follows that $Y$ is a $(g+md_{\Gamma })$-supermartingale on $[0,T]$.
Conversely, if for each $m\geq 0$, $Y$ is a $(g+md_{\Gamma })$%
-supermartingale on $[0,T]$ i.e. $\mathcal{E}_{s,t}^{g+md_{\Gamma
}}[Y_{t}]=y_{s}^{t,m}\leq Y_{s}$. When we let $m\rightarrow \infty $, by the
monotonic limit theorem in \cite{P99}, $y_{\cdot }^{t,m}$ converges to $%
\mathcal{E}_{\cdot ,t}^{g_{\Gamma }}[Y_{t}]$, which is the $g_{\Gamma }$%
--solution. We thus have $\mathcal{E}_{s,t}^{g_{\Gamma }}[Y_{t}]\leq Y_{s}$,
on $[0,T]$. This implies that $Y$ is a $g_{\Gamma }$-supermartingale.\hfill $%
\Box $\medskip

We have the following $g_{\Gamma }$-supermartingale decomposition theorem.

\begin{theorem}
\label{decomp1}Let $Y$ be a right continuous $g_{\Gamma }$-supermartingale
on $[0,T]$. Then there exists a unique RCLL increasing process $A\in \mathbf{%
A}_{\mathcal{F}}^{2}(0,T)$, such that $Y$ is a $g_{\Gamma }$--supersolution,
namely,
\begin{eqnarray*}
y_{t}
&=&Y_{T}+\int_{t}^{T}g(s,y_{s},z_{s})ds+A_{T}-A_{t}-\int_{t}^{T}z_{s}dB_{s},
\\
&&\
\begin{array}{c}
d_{\Gamma _{t}}(y_{t},z_{t})=0\mbox{, a.e. a.s..}
\end{array}
\end{eqnarray*}
\end{theorem}

\noindent\textbf{Proof. } For each fixed $m\geq 0$, we consider the
solution $(y^{m},z^{m},A^{m})\in
\mathbf{D}_{\mathcal{F}}^{2}(0,T)\times \mathbf{L}_{\mathcal{F}}^{2}(0,T;%
\mathbb{R}^{d})\times \mathbf{A}_{\mathcal{F}}^{2}(0,T)$ of the following
reflected BSDE, with the $g_{\Gamma }$--supermartingale as the lower
obstacle:
\begin{eqnarray}
y_{t}^{m}
&=&Y_{T}+\int_{t}^{T}g^{m}(s,y_{s}^{m},z_{s}^{m})ds+A_{T}^{m}-A_{t}^{m}-%
\int_{t}^{T}z_{s}^{m}dB_{s},  \label{pesupm} \\
dA^{m} &\geq &0,\;\;y_{t}^{m}\leq
Y_{t},\;\;\int_{0}^{T}(Y_{t}-y_{t}^{m})dA_{t}^{m}=0.  \nonumber
\end{eqnarray}
where $g^{m}(t,y,z):=(g+md_{\Gamma })(t,y,z)$. By Proposition \ref{gmphi},
this $g_{\Gamma }$--supermartingale is also a $g^{m}$--supermartingale for
each $m$. It follows from the $g$--supermartigale decomposition theorem (see
\cite{P99}) that $y_{t}^{m}\equiv Y_{t}$. Thus $z^{m}$ is invariant in $m$%
: $z^{m}\equiv Z\in \mathbf{L}_{\mathcal{F}}^{2}(0,T;\mathbb{R}^{d})$ and
the above equation (\ref{pesupm}) can be written
\[
Y_{t}=Y_{T}+\int_{t}^{T}(g+md_{\Gamma
})(s,Y_{s},Z_{s})ds+A_{T}^{m}-A_{t}^{m}-\int_{t}^{T}Z_{s}dB_{s}.
\]
Consequently, for all $m\geq 0$, notice that $A^{m}$ is a positive process,
we have
\begin{eqnarray*}
0 &\leq &m\int_{0}^{T}d_{\Gamma }(s,Y_{s},Z_{s})ds \\
&\leq &\left[
Y_{0}-Y_{T}+\int_{0}^{T}Z_{s}dB_{s}-\int_{0}^{T}g(s,Y_{s},Z_{s})ds\right]
^{+}\in \mathbf{L}^{2}(\mathcal{F}_{T}).
\end{eqnarray*}
From this it follows immediately $\int_{0}^{T}d_{\Gamma }(s,Y_{s},Z_{s})ds=0$%
. Thus $A^{m}$ is also invariant in $m$: $A^{m}=A\in \mathbf{D}_{\mathcal{F%
}}^{2}(0,T)$. We thus complete the proof. \hfill $\Box $\medskip

\subsection{$g_{\Gamma }$--submartingale decomposition theorem}

We now consider the decomposition theorem of a given $g_{\Gamma }$%
--submartingale $Y\in \mathbf{D}_{\mathcal{F}}^{2}(0,T)$. In
\cite{PX2006}, we have proved a $g_{\Gamma }$-submartingale
decomposition theorem under assumptions of $Y_{t-}\geq Y_{t}$ with

\noindent \textbf{(H)} There exists a quadruple $(y^{*},z^{*},A^{*},K^{*})%
\in \mathbf{D}_{\mathcal{F}}^{2}(0,T)\times \mathbf{L}_{\mathcal{F}%
}^{2}(0,T)\times (\mathbf{A}_{\mathcal{F}}^{2}(0,T))^2$, satisfying
\begin{eqnarray*}
y_{t}^{*}
&=&Y_{T}+%
\int_{t}^{T}g(s,y_{s}^{*},z_{s}^{*})ds+(A_{T}^{*}-A_{t}^{*})-(K_{T}^{*}-K_{t}^{*})-\int_{t}^{T}z_{s}^{*}dB_{s},
\\
&&
\begin{array}{c}
d_{\Gamma _{t}}(y_{t}^{*},z_{t}^{*})=0\mbox{, a.s. a.e.}
\end{array}
\\
y_{t}^{*} &\leq &Y_{t}\mbox{, \ \ }%
\int_{0}^{T}(y_{t-}^{*}-Y_{t-})dK_{t}^{*}=0\mbox{, a.s..}
\end{eqnarray*}

\begin{remark}
A necessary condition for (\textbf{H) }holding is $\Gamma _{t}\cap (-\infty
,Y_{t}]\times \mathbb{R}^{d}\neq \emptyset $.
\end{remark}

Here we partly generalize this result and try to get rid of
assumption \textbf{(H)}.

\begin{theorem}
\label{Sub-dec}Assume $\Gamma _{t}$ only depends on $z$. Let $Y\in \mathbf{D}%
_{\mathcal{F}}^{2}(0,T)$ be a $g_{\Gamma }$-submartingale on $[0,T]$ and for
stopping times $\sigma $, $\tau $ valued in $[0,T]$, with $\sigma \leq \tau $%
, such that
\begin{equation}
\mathcal{E}_{\sigma ,\tau }^{g_{\Gamma }}[Y_{\tau }-(Y_{\tau }-Y_{\tau
-})^{+}]\geq Y_{\sigma }.  \label{cond-sub}
\end{equation}
Then there exists a unique continuous increasing process $K\in \mathbf{A}_{\mathcal{F}}^{2}(0,T)$, such that the triple $(Y-K,Z,A)\in \mathbf{D}_{%
\mathcal{F}}^{2}(0,T)\times \mathbf{L}_{\mathcal{F}}^{2}(0,T;\mathbb{R}%
^{d})\times \mathbf{A}_{\mathcal{F}}^{2}(0,T)$ is the $g^{K}{}_{\Gamma }$%
--solution with terminal condition $Y_{T}-K_{T}$, i.e. for $t\in [0,T],$%
\begin{eqnarray*}
Y_{t}-K_{t}
&=&Y_{T}-K_{T}+\int_{t}^{T}g^{K}(s,Y_{s}-K_{s},Z_{s})ds+(A_{T}-A_{t})-%
\int_{t}^{T}Z_{s}dB_{s}, \\
Z_{t} &\in &\Gamma _{t},\;\;\;dP\times dt\mbox{-a.s.}
\end{eqnarray*}
where
\[
g^{K}(t,y,z):=g(t,y+K_{t},z),\;(t,y,z)\in [0,T]\times \mathbb{R}\times %
\mathbb{R}^{d}.
\]
\end{theorem}

\noindent\textbf{Proof. } Consider the BSDE$(Y_{T},g_{\Gamma })$
with reflecting upper obstacle $Y$.
From Theorem \ref{existru}, we know that there exists a quadruple $%
(y,Z,A,K)\in \mathbf{D}_{\mathcal{F}}^{2}(0,T)\times \mathbf{L}_{\mathcal{F}%
}^{2}(0,T;\mathbb{R}^{d})\times (\mathbf{A}_{\mathcal{F}}^{2}(0,T))^{2}$
\begin{eqnarray}
y_{t}
&=&Y_{T}+\int_{t}^{T}g(s,y_{s},Z_{s})ds+A_{T}-A_{t}-(K_{T}-K_{t})-%
\int_{t}^{T}Z_{s}dB_{s},  \label{RCBSDEx} \\
Z_{t} &\in &\Gamma _{t}\mbox{,\ }dP\times dt\mbox{-a.s.,\ }dA\geq
0,\;\;dK\geq
0,\;\;\mathcal{V}_{[0,t]}[A-K]=\mathcal{V}_{[0,t]}[A+K],
\nonumber \\
y_{t} &\leq
&Y_{t},\;\;\int_{0}^{T}(y_{s-}-Y_{s-})dK_{s}=0,\;\mbox{a.s..}
\nonumber
\end{eqnarray}

We want to prove that $y\equiv Y$. It is sufficient to prove
$y_{t}\geq Y_{t} $. For each $\delta >0$, we define stopping times
\begin{eqnarray*}
\sigma ^{\delta } &:&=\inf \{t,y_{t}\leq Y_{t}-\delta \}\wedge T, \\
\tau &:&=\inf \{t\geq \sigma ^{\delta }:\;y_{t}\geq Y_{t}\}.
\end{eqnarray*}
If $P(\sigma ^{\delta }<T)=0$ for all $\delta >0$, the proof is done; if it
is not such case, there exists a $\delta >0$, such that $P(\sigma ^{\delta
}<T)>0$. So we have $\sigma ^{\delta }<\tau \leq T$. Since $y$ and $Y$ are
RCLL, $y_{\sigma ^{\delta }}\leq Y_{\sigma ^{\delta }}-\delta $ and $y_{\tau
}\geq Y_{\tau }$. So $y_{\tau }=Y_{\tau }$. By the integral equality in (\ref
{RCBSDEx}), we get $K_{\tau -}=K_{\sigma ^{\delta }}$.

Since $\mathcal{V}_{[0,t]}[A-K]=\mathcal{V}_{[0,t]}[A+K]$,
$\bigtriangleup A_{\tau }\cdot \bigtriangleup K_{\tau }=0$. From the
integral equality in (\ref
{RCBSDEx}), we know that $(y_{s-}-Y_{s-})(K_{s}-K_{s-})=0$. So if $%
\bigtriangleup K_{\tau }\neq 0$, then $y_{\tau -}=Y_{\tau -}$ and $%
\bigtriangleup A_{\tau }=0$, which implies $\bigtriangleup K_{\tau
}=(y_{\tau }-y_{\tau -})^{+}=(Y_{\tau }-Y_{\tau -})^{+}$.

Define
\begin{eqnarray*}
\overline{Y}_{t} &=&Y_{t}1_{[0,\tau )}(t)+(Y_{\tau }-(Y_{\tau }-Y_{\tau
-})^{+})1_{[\tau ,T]}(t), \\
\overline{g}(t,y,z) &=&g(t,y,z)1_{[0,\tau )}(t), \\
\overline{\Gamma }_{t} &=&\Gamma _{t}1_{[0,\tau ]}+R^{1\times d}1_{[\tau
,T]}(t).
\end{eqnarray*}
Then for $0\leq s\leq t\leq T$, $\mathcal{E}_{s,t}^{\overline{g},\overline{%
\Gamma }}[\overline{Y}_{t}]=\mathcal{E}_{s\wedge \tau ,t\wedge \tau
}^{g,\Gamma }[Y_{t\wedge \tau }-(Y_{\tau }-Y_{\tau -})^{+}1_{[\tau ,T]}(t)]$.

Consider two stopping times $0\leq \sigma _{1}\leq \sigma _{2}\leq
T$, then with (\ref{cond-sub}) we have
\begin{eqnarray*}
\mathcal{E}_{\sigma _{1},\sigma _{2}}^{\overline{g}_{\overline{\Gamma }}}[%
\overline{Y}_{\sigma _{2}}] &=&\mathcal{E}_{\sigma _{1}\wedge \tau ,\sigma
_{2}\wedge \tau }^{g,\Gamma }[Y_{\sigma _{2}\wedge \tau }-(Y_{\tau }-Y_{\tau
-})^{+}1_{[\tau ,T]}(\sigma _{2})] \\
&\geq &\mathcal{E}_{\sigma _{1}\wedge \tau ,\sigma _{2}\wedge \tau
}^{g,\Gamma }[Y_{\sigma _{2}\wedge \tau }-(Y_{\sigma _{2}\wedge \tau
}-Y_{\sigma _{2}\wedge \tau -})^{+}] \\
&\geq &Y_{\sigma _{1}\wedge \tau }\geq \overline{Y}_{\sigma _{1}}.
\end{eqnarray*}
So $\overline{Y}$ is a $\overline{g}_{\overline{\Gamma
}}$-submartingale. Define
\begin{eqnarray*}
\overline{y}_{t} &=&y_{t}1_{[0,\tau )}(t)+(y_{\tau }-(y_{\tau }-y_{\tau
-})^{+})1_{[\tau ,T]}(t), \\
\overline{K}_{t} &=&K_{t}1_{[0,\tau )}(t)+K_{\tau -}1_{[\tau ,T]}(t), \\
\overline{A}_{t} &=&A_{t}1_{[0,\tau ]}(t)+A_{\tau}1_{(\tau ,T]}(t),
\;\; \overline{Z}_{t} =Z_{t}1_{[0,\tau ]}(t),
\end{eqnarray*}
Then for $0\leq t\leq T$, we have
\begin{eqnarray*}
\overline{y}_{t} &=&\overline{Y}_{T}+\int_{t}^{T}g(s,\overline{y}_{s},%
\overline{z}_{s})ds+(\overline{A}_{T}-\overline{A}_{t})-(\overline{K}_{T}-%
\overline{K}_{t})-\int_{t}^{T}\overline{Z}_{s}dB_{s}, \\
\overline{Z}_{t} &\in &\overline{\Gamma }_{t}\mbox{,\ }dP\times dt\mbox{%
-a.s.,\ }\ \ d\overline{A}\geq 0,\;\;d\overline{K}\geq 0,\;\;\mathcal{V}%
_{[0,t]}[\overline{A}-\overline{K}]=\mathcal{V}_{[0,t]}[\overline{A}+%
\overline{K}], \\
\overline{y}_{t} &\leq &\overline{Y}_{t},\;\;\;\int_{0}^{T}(\overline{y}%
_{s-}-\overline{Y}_{s-})d\overline{K}_{s}=0.
\end{eqnarray*}
Notice that $\overline{K}_{\tau }=K_{\tau -}=K_{\sigma ^{\delta }}=\overline{%
K}_{\sigma ^{\delta }}$, with $\overline{y}_{\tau }=y_{\tau }-\bigtriangleup
K_{\tau }=Y_{\tau }-\bigtriangleup K_{\tau }=\overline{Y}_{\tau }$, we get
\begin{eqnarray*}
\overline{y}_{\sigma ^{\delta }} &=&\overline{Y}_{\tau }+\int_{\sigma
^{\delta }}^{\tau }g(s,\overline{y}_{s},\overline{z}_{s})ds+(\overline{A}%
_{\tau }-\overline{A}_{\sigma ^{\delta }})-\int_{\sigma ^{\delta }}^{\tau }%
\overline{Z}_{s}dB_{s}, \\
(\overline{y}_{t},\overline{Z}_{t}) &\in &\overline{\Gamma }_{t},\ \ d%
\overline{A}\geq 0.
\end{eqnarray*}
Since $y+K$ is the $g_{\Gamma }$-solution of constraint BSDE $(Y_{\tau
},g,\Gamma )$ on $[\sigma ^{\delta },\tau ]$, $\overline{y}$ is the $%
g_{\Gamma }$-solution of constraint BSDE $(\overline{Y}_{\tau },\overline{g},%
\overline{\Gamma })$ on $[\sigma ^{\delta },\tau ]$. So with the fact that $%
\overline{Y}$ is a $\overline{g}_{\overline{\Gamma }}$-submartingale, we get
\[
\overline{y}_{\sigma ^{\delta }}=\mathcal{E}_{\sigma ^{\delta },\tau }^{%
\overline{g}_{\overline{\Gamma }}}[\overline{y}_{\tau }]=\mathcal{E}_{\sigma
^{\delta },\tau }^{\overline{g}_{\overline{\Gamma }}}[\overline{Y}_{\tau
}]\geq \overline{Y}_{\sigma ^{\delta }}.
\]
But $\overline{y}_{\sigma ^{\delta }}=y_{\sigma ^{\delta }}$, $\overline{Y}%
_{\sigma ^{\delta }}=Y_{\sigma ^{\delta }}$, this introduces a
contradiction. $\square $

\begin{remark}
This result can easy cover the decomposition theorem in
\cite{PX2006} when $\Gamma$ does not depend on $y$. In
fact, from the condition $Y_{t-}\geq Y_{t}$, we know that $%
(Y_{t}-Y_{t-})^{+}=0$, so $\mathcal{E}_{s,t}^{g_{\Gamma
}}[Y_{t}-(Y_{t}-Y_{t-})^{+}]=\mathcal{E}_{s,t}^{g_{\Gamma }}[Y_{t}]\geq Y_{s}
$.
\end{remark}

\begin{remark}
We can prove the same result for general case when $\Gamma $ also
depends on $y$ under the assumption \textbf{(H)}. This assumption is
not easy to verify. However it is required by the existence of
$g_{\Gamma }$-reflected solution associated to $(Y_{T},g,\Gamma )$
with an upper barrier $Y$, when the constraint $\Gamma$ depends on
$y$.
\end{remark}

Although in Theorem \ref{Sub-dec}, we remove assumption
\textbf{(H)}, sometimes the assumption (\ref{cond-sub}) in is not
easy either. In the following result, we do not need to assume
(\ref{cond-sub}), but we need more assumptions on $g$.

\begin{theorem}
\label{sub-dec2}Let $Y\in \mathbf{D}_{\mathcal{F}}^{2}(0,T)$ be a
$g_{\Gamma }$-submartingale on $[0,T]$. Suppose $g$ and $\Gamma $ do
not depend on $y$, $g(t,0)=0$ and $0\in \Gamma _{t}$. Then there
exists a unique continuous increasing process $K$ with
$E[K_{T}^{2}]<\infty $, such that the same decomposition result of
theorem \ref{Sub-dec} holds.
\end{theorem}

\noindent\textbf{Proof. }
As in the proof of theorem \ref{Sub-dec}, we consider the BSDE$%
(Y_{T},g_{\Gamma })$ with reflecting upper obstacle $Y$. From Theorem \ref
{existru}, we know that there exists a quadruple $(y,Z,A,K)\in \mathbf{D}_{%
\mathcal{F}}^{2}(0,T)\times \mathbf{L}_{\mathcal{F}}^{2}(0,T;\mathbb{R}%
^{d})\times (\mathbf{A}_{\mathcal{F}}^{2}(0,T))^{2}$ such that
\begin{eqnarray}
y_{t}
&=&Y_{T}+\int_{t}^{T}g(s,Z_{s})ds+A_{T}-A_{t}-(K_{T}-K_{t})-%
\int_{t}^{T}Z_{s}dB_{s},  \label{RCBSDE2} \\
Z_{t} &\in &\Gamma _{t}\mbox{,\ }dP\times dt\mbox{-a.s.,\ }\ \ \mbox{\ \ }%
dA\geq 0,\;\;dK\geq 0,\;\;\mathcal{V}_{[0,t]}[A-K]=\mathcal{V}_{[0,t]}[A+K],
\nonumber \\
y_{t} &\leq
&Y_{t},\;\;\int_{0}^{T}(y_{s-}-Y_{s-})dK_{s}=0,\;\mbox{a.s..}
\nonumber
\end{eqnarray}

As before we want to prove that $y\equiv Y$. It suffices to prove
$y_{t}\geq Y_{t} $. For each $\delta >0$, define stopping times
\begin{eqnarray*}
\sigma ^{\delta } &:&=\inf \{t,y_{t}\leq Y_{t}-\delta \}\wedge T, \\
\tau ^{\delta } &:&=\inf \{t\geq \sigma ^{\delta }:\;y_{t}\geq Y_{t}-\frac{%
\delta }{2}\}.
\end{eqnarray*}
If $P(\sigma ^{\delta }<T)=0$ for all $\delta >0$, the proof is
done; if it is not, there exists a $\delta >0$, such that $P(\sigma
^{\delta }<T)>0$. So we have $\sigma ^{\delta }<\tau ^{\delta }\leq
T$. Since $y$ and
$Y$ are RCLL, $y_{\sigma ^{\delta }}\leq Y_{\sigma ^{\delta }}-\delta $ and $%
y_{\tau ^{\delta }}\geq Y_{\tau ^{\delta }}-\frac{\delta }{2}$. By the
integral equality in (\ref{RCBSDE2}), we get $K_{\tau ^{\delta }}=K_{\sigma
^{\delta }}$.

Thanks to proposition \ref{tr-con}-(i), we know that
$\mathcal{E}^{g_{\Gamma }}[\cdot ]$ has translation invariant
property. So $\mathcal{E}_{\sigma ^{\delta
},\tau ^{\delta }}^{g_{\Gamma }}[Y_{\tau ^{\delta }}-\frac{\delta }{2}]=%
\mathcal{E}_{\sigma ^{\delta },\tau ^{\delta }}^{g_{\Gamma
}}[Y_{\tau ^{\delta }}]-\frac{\delta }{2}$.

While on the interval $[\sigma ^{\delta },\tau ^{\delta }]$,
\begin{eqnarray*}
y_{\sigma ^{\delta }} &=&y_{\tau ^{\delta }}+\int_{\sigma ^{\delta }}^{\tau
^{\delta }}g(s,Z_{s})ds+A_{\tau ^{\delta }}-A_{\sigma ^{\delta
}}-\int_{\sigma ^{\delta }}^{\tau ^{\delta }}Z_{s}dB_{s}, \\
Z_{t} &\in &\Gamma _{t}\mbox{,\ }dP\times dt\mbox{-a.s.,\ \ \
}dA\geq 0.
\end{eqnarray*}
So we have
\[
y_{\sigma ^{\delta }}=\mathcal{E}_{\sigma ^{\delta },\tau ^{\delta
}}^{g_{\Gamma }}[y_{\tau ^{\delta }}]\geq \mathcal{E}_{\sigma ^{\delta
},\tau ^{\delta }}^{g_{\Gamma }}[Y_{\tau ^{\delta }}-\frac{\delta }{2}]=%
\mathcal{E}_{\sigma ^{\delta },\tau ^{\delta }}^{g_{\Gamma }}[Y_{\tau
^{\delta }}]-\frac{\delta }{2}\geq Y_{\sigma ^{\delta }}-\frac{\delta }{2}.
\]
This introduces a contradiction. So result follows. $\square $

\begin{remark}
\label{sub-dec3}By proposition \ref{tr-con}-(ii), we can prove the same
results, for the case when $g(t,y,z)=g_{1}(t,z)+ay$ with $g_{1}(t,z)$ is bounded, and $%
\Gamma $ only depends on $z$, with $0\in \Gamma _{t}$.
\end{remark}

\section{Applications of $g_{\Gamma }$-reflected BSDEs: American option
pricing in incomplete market}

We follow the idea of El Karoui et al.(1997, \cite{EPQ}). Consider
the strategy wealth portfolio $(Y_{t},\pi _{t})$ as a pair of
adapted processes
in $\mathbf{L}_{\mathcal{F}}^{2}(0,T)\times \mathbf{L}_{\mathcal{F}}^{2}(0,T;%
\mathbb{R}^{d})$ which satisfy the following BSDE
\[
-dY_{t}=g(t,Y_{t},\pi _{t})dt-\pi _{t}^{\tau}\sigma _{t}dB_{t},
\]
where $g$ is $\mathbb{R}$-valued, convex with respect to $(y,\pi )$,
and satisfy Lipschitz condition (\ref{Lip}). We suppose that the
volatility matrix $\sigma $ of $n$ risky assets is invertible and
 $(\sigma _{t})^{-1}$ is bounded.

In complete market, we are concerned with the problem of pricing an American
contingent claim at each time $t$, which consists of the selection of a
stopping time $\tau \in \mathcal{T}_{t}$ (the set of stopping times valued
in $[t,T]$) and a payoff $S_{\tau }$ on exercise if $\tau <T$ and $\xi $ if $%
\tau =T$. Here $(S_{t})$ is a continuous process satisfying $%
E[\sup_{t}(S_{t}^{+})^{2}]<\infty $. Set
\[
\widetilde{S}_{s}=\xi 1_{\{s=T\}}+S_{s}1_{\{s<T\}}.
\]
Then the price of the American contingent claim $(\widetilde{S}_{s},0\leq
s\leq T)$ at time $t$ is given by
\[
Y_{t}=ess\sup_{\tau \in \mathcal{T}_{t}}Y_{t}(\tau ,\widetilde{S}_{\tau }).
\]
Moreover the price $(Y_{t},0\leq t\leq T)$ corresponds to the unique
solution of the reflected BSDE associated with terminal condition $\xi $, generator $%
g $ and obstacle $S$, i.e. there exists $(\pi _{t})\in \mathbf{L}_{\mathcal{F%
}}^{2}(0,T;\mathbb{R}^{d})$ and $(A_{t})$ an increasing continuous
process with $A_{0}=0$ such that
\begin{eqnarray}
-dY_{t} &=&g(s,Y_{t},\pi _{t})ds+dA_{t}-\pi _{t}^{\tau}\sigma
_{t}dB_{t},Y_{T}=\xi ,  \nonumber \\
Y_{t}\geq S_{t} &,&0\leq t\leq T,\;\;\;\int_{0}^{T}(Y_{t}-S_{t})dA_{t}=0.
\nonumber
\end{eqnarray}
Furthermore, the stopping time $D_{t}=\inf (t\leq s\leq T\mid
Y_{s}=S_{s})\wedge T$ is optimal, that is
\[
Y_{t}=Y_{t}(D_{t},\widetilde{S}_{D_{t}}).
\]

Now we consider in the incomplete market, i.e. there is a constraint
on portfolio $\pi _{t}\in \Gamma _{t}$, where $\Gamma _{t}$ is a
closed subset
of $\mathbb{R}^{d}$, how to price the American contingent claim $(\widetilde{%
S}_{s},0\leq s\leq T)$. Lucky, with the results in former sections, we have
the following results:

\begin{theorem}
If $\xi $ is attainable, i.e. there exists a couple $(Y^{\prime
},\pi ^{\prime })$ with $\pi _{t}^{\prime }\in \Gamma _{t}$,
$t$-a.e. which replicate $\xi $, then the price process $Y$ of
American option in the incomplete market is the $g_{\Gamma
}$--solution reflected by the lower
obstacle $L$, i.e. there exist a process $\pi _{t}\in \Gamma _{t}$, $dP\times dt$%
-a.s., and increasing continuous processes $A$ and $\bar{A}$, such
that
\begin{eqnarray}
Y_{t} &=&\xi +\int_{t}^{T}g(s,Y_{s},\pi _{s})ds+A_{T}-A_{t}+\overline{A}_{T}-%
\overline{A}_{t}-\int_{t}^{T}\pi _{s}^{\tau}\sigma _{s}dB_{s},  \label{CAOP} \\
Y_{t} &\geq &S_{t},0\leq t\leq T,\;\;\;\;\int_{0}^{T}(Y_{t}-S_{t})d\overline{%
A}_{t}=0.\   \nonumber
\end{eqnarray}
Furthermore, the stopping time $D_{t}=\inf (t\leq s\leq T\mid
Y_{s}=S_{s})\wedge T$ is still optimal.
\end{theorem}

\textit{Sketch of the proof. }Thanks to the results of \cite{EPQ}
and \cite {W02}, we know that the method of auxiliary market in
\cite{CK92} and \cite {CK93} is equivalent to the penalization
equations associated to $(\xi ,f+nd_{\Gamma _{t}},S)$, then let
$n\rightarrow \infty $, we may get the price.
By theorem \ref{existrl}, since $\xi $ is attainable, the result follows. $%
\square $

\subsection{Some examples of American call option}

We study the American call option, set $S_{t}=(X_{t}-k)^{+}$, $\xi
=(X_{T}-k)^{+}$, where $X$ is the price of underlying stock and $k$ is the
strike price. More precisely, $X$ is the solution of
\begin{equation}
X_{t}=x_{0}+\int_{0}^{t}\mu _{s}X_{s}ds+\int_{0}^{t}\sigma _{s}X_{s}dB_{s}.
\label{stprice}
\end{equation}
Correspondingly, in (\ref{CAOP}) $g$ is a linear function
\[
g(t,y,\pi )=-r_{t}y-(\mu _{t}-r_{t})\pi ^{\tau}\sigma _{t}.
\]

\begin{proposition}
\label{AOPcall}If $\xi $ is attainable, then the maturity time of American
call option in incomplete market is still $T$ .
\end{proposition}

\noindent\textbf{Proof. }
We have that $Y_{t}\geq Y_{t}^{0}$, $\overline{A}_{t}\leq \overline{A}%
_{t}^{0}$, $t\in [0,T]$, comparing (\ref{CAOP}) and $Y^{0}$, where
$Y^{0}$ is the price process of American call option without
constraint, which satisfies a reflected BSDE
\begin{eqnarray*}
Y_{t}^{0} &=&\xi +\int_{t}^{T}g(s,Y_{s}^{0},\pi _{s}^{0})ds+\overline{A}%
_{T}^{0}-\overline{A}_{t}^{0}-\int_{t}^{T}(\pi _{s}^{0})^{\tau}\sigma _{s}dB_{s}, \\
Y_{t}^{0} &\geq &S_{t}\mbox{, \ \ \ }\int_{0}^{T}(Y_{t}^{0}-S_{t})d\overline{%
A}_{t}^{0}=0.\;
\end{eqnarray*}

Since American call option always exercises at terminal time $T$,
which implies $\overline{A}_{t}^{0}=0$ and $D_{t}^{0}=T$, where
$D_{t}^{0}=\inf (t\leq s\leq T\mid Y_{s}^{0}=S_{s})\wedge T$. So we
have $Y_{t}^{0}>S_{t}$
on $[0,T)$. It follows that $Y_{t}\geq Y_{t}^{0}>S_{t}$ on $[0,T)$ and $%
\overline{A}_{t}\leq \overline{A}_{t}^{0}=0$, $t\in [0,T]$. Then $D_{t}=T$. $%
\square $

From this proposition, we know that there is no difference between the
American call option and European call option even in incomplete market.

\begin{example}
No short-selling: In this case $\Gamma _{t}=[0,\infty )$, for $t\in [0,T]$.
We set $d=1$. By the proposition \ref{AOPcall} and Example 7.1 in \cite{CK93}%
, the price process of the American call option takes same value as
European call option. This means that the constraint $K=[0,\infty )$
does not make any difference.
\end{example}

In fact, we have a more general result.

\begin{proposition}
\label{monoCall}Consider the constraint $\Gamma _{t}=[0,\infty )$, for $t\in
[0,T]$. If $\xi =\Phi (X_{T})$, $S_{t}=l(X_{t})$, where $\Phi $, $l:%
\mathbb{R}\mathbf{\rightarrow }\mathbb{R}$ are both increasing in
$x$, and $\sigma $ satisfies the uniformly elliptic condition, then
the price
process $Y$ takes same value as in complete market, i.e. the constraint $%
\Gamma $ does not influence the price.
\end{proposition}

\noindent\textbf{Proof. }
It is sufficient to prove that $\overline{\pi }_{t}\geq 0$, where $(%
\overline{Y},\overline{\pi },\overline{A})$ is the solution of following
reflected BSDE
\begin{eqnarray}
\overline{Y}_{t} &=&\Phi (X_{T})+\int_{t}^{T}g(s,\overline{Y}_{s},\overline{%
\pi }_{s})ds+\overline{A}_{T}-\overline{A}_{t}-\int_{t}^{T}\overline{\pi }%
_{s}^{\tau}\sigma _{s}dB_{s},  \label{REcall} \\
\overline{Y}_{t} &\geq &l(X_{t})\mbox{, \ \ \ }\int_{0}^{T}(\overline{Y}%
_{t}-l(X_{t}))d\overline{A}_{t}=0.\;  \nonumber
\end{eqnarray}

We put $(X_{s}^{t,x},\overline{Y}_{s}^{t,x},\overline{\pi
}_{s}^{t,x},\overline{A}_{s}^{t,x})_{t\leq s\leq T}$ under Markovian
framework. Define
\[
u(t,x)=\overline{Y}_{t}^{t,x},
\]
then by \cite{EKPPQ}, we know that $u$ is the viscosity solution of the PDE
with an obstacle $l$,
\begin{eqnarray*}
\min \{u(t,x)-l(x),-\frac{\partial u}{\partial t}-\mathcal{L}%
u-g(t,x,u,\nabla u\sigma )\} &=&0, \\
u(T,x) &=&\Phi (x),
\end{eqnarray*}
where $\mathcal{L}=\frac{1}{2}(\sigma _{s})^{2}\frac{\partial ^{2}}{\partial
x\partial x}+\mu \frac{\partial }{\partial x}$. Since $(\overline{\pi }%
_{r}^{t,x})^{\tau}\sigma _{r}=\nabla u\sigma (r,X_{r}^{t,x})$, and
$\sigma$ is uniformly elliptic, we only need to prove that $\nabla
u(t,x)$ is non-negative. Indeed, it is easy to obtain by comparison
theorem.
For $x_{1}$, $x_{2}\in \mathbb{R}$, with $x_{1}\geq x_{2}$, $%
X_{s}^{t,x_{1}}\geq X_{s}^{t,x_{2}}$. It follows that $\Phi
(X_{T}^{t,x_{1}})\geq \Phi (X_{T}^{t,x_{2}})$ and
$l(X_{s}^{t,x_{1}})\geq
l(X_{s}^{t,x_{2}})$ in view of assumptions. By comparison theorem of BSDE, $%
\overline{Y}_{t}^{t,x_{1}}\geq \overline{Y}_{t}^{t,x_{2}}$, which implies $%
u(t,x_{1})\geq u(t,x_{2})$. So $\nabla u(t,x)\geq 0$, it follows that $%
\overline{\pi }_{t}^{t,x}\geq 0$. $\square $

\subsection{Some examples of American put option}

In this case, we set $S_{t}=(k-X_{t})^{+}$, $\xi =(k-X_{T})^{+}$, where $X$
is the price of underlying stock as in (\ref{stprice}) and $k$ is the strike
price. Similarly to proposition \ref{monoCall}, we have

\begin{proposition}
\label{monoPut}Consider the constraint $\Gamma _{t}=(-\infty ,0]$, for $t\in
[0,T].$ If $\xi =\Phi (X_{T})$, $S_{t}=l(X_{t})$, where $\Phi $, $l:%
\mathbb{R}\mathbf{\rightarrow }\mathbb{R}$ are both decreasing
functions, and $\sigma $ satisfies uniformly elliptic condition,
then the price process $Y$ takes same value as in complete market,
i.e. the constraint $\Gamma $ has no influence on price process.
\end{proposition}

\noindent\textbf{Proof. } Similar to the proof of proposition
\ref{monoCall}, it is sufficient to prove that $\overline{\pi
}_{t}\leq 0$. With the helps of viscosity solution, we get the
result. $\square $

\begin{example}
No borrowing: $\Gamma _{t}=(-\infty ,Y_{t}]$. Obviously, $Y_{t}\geq
0$, in view of $Y_{t}\geq S_{t}\geq 0$. So $\Gamma _{t}\supset
(-\infty ,0]$, by proposition \ref{monoPut}, we know that the price
process $Y$ takes same value as in complete market. This means that
to replicate an American put option, we don't need to borrow money.
\end{example}

\begin{example}
No short-selling: $\Gamma _{t}=[0,\infty )$, for $t\in [0,T]$. Then the
pricing process $Y$ with hedging $\pi $ satisfying
\begin{eqnarray*}
Y_{t} &=&\xi +\int_{t}^{T}g(s,Y_{s},\pi _{s})ds+A_{T}-A_{t}+\overline{A}_{T}-%
\overline{A}_{t}-\int_{t}^{T}\pi _{s}^{*}\sigma _{s}dB_{s}, \\
Y_{t} &\geq &S_{t},0\leq t\leq T,\;\int_{0}^{T}(Y_{t}-S_{t})d\overline{A}%
_{t}=0,\ \pi _{t}\geq 0,t\mbox{-a.e..}
\end{eqnarray*}
Notice that $S_{t}=(k-X_{t})^{+}<k$. So the $g_{\Gamma }$-solution of the
above equation is
\begin{eqnarray*}
Y_{t} &=&\left\{
\begin{array}{ll}
k, & t\in [0,T) \\
(k-X_{T})^{+}, & t=T
\end{array}
;\right.  \\
\pi _{t} &=&0, \\
A_{t} &=&\left\{
\begin{array}{ll}
k\int_{0}^{t}r_{s}ds, & t\in [0,T) \\
k\int_{0}^{T}r_{s}ds+k-(k-X_{T})^{+}, & t=T
\end{array}
;\right.  \\
\overline{A}_{t} &=&0.
\end{eqnarray*}
In particular, $Y_{0}=k$, which is the price of American put option under
'no short-selling' constraint.
\end{example}

\section{Appendix}

In appendix, we recall some results of $g_{\Gamma}$-solution in
\cite{P99}, and proves some comparison results of
$g_{\Gamma}$-solution. In \cite{P99}, $\Gamma $ is defined as
\[
\Gamma _{t}(\omega )=\{(y,z)\in \mathbb{R}^{1+d}:\Phi (\omega ,t,y,z)=0\}.
\]
where $\Phi $ is a nonnegative, measurable Lipschitz function and $\Phi
(\cdot ,y,z)\in \mathbf{L}_{\mathcal{F}}^{2}(0,T)$, for $(y,z)\in \mathbb{%
R\times R}^{d}$. Under the following assumption, the existence of the
smallest solution is proved.

The following theorem of the existence of the smallest solution was obtained
in \cite{P99}.

\begin{theorem}\label{exist}
Suppose that the function $g$ satisfies (\ref{Lip}) and the constraint $%
\Gamma $ satisfies (\ref{Gamma}). We assume that there is at least one $%
\Gamma $--constrained $g$--supersolution $y^{\prime }\in \mathbf{D}_{%
\mathcal{F}}^{2}(0,T)$:
\begin{eqnarray}
y_{t}^{\prime } &=&X^{\prime }+\int_{t}^{T}g(s,y_{s}^{\prime
},z_{s})ds+A_{T}^{\prime }-A_{t}^{\prime }-\int_{t}^{T}z_{s}^{\prime }dB_{s},
\label{no-empty1} \\
A &\in &\mathbf{A}_{\mathcal{F}}^{2}(0,T)\mbox{ ,\ }(y_{t}^{\prime
},z_{t}^{\prime })\in \Gamma _{t},\;t\in [0,T]\mbox{, a.s. a.e.}
\nonumber
\end{eqnarray}
Then, for each $X\in \mathbf{L}^{2}(\mathcal{F}_{T})$ with $X\leq
X^{\prime } $, a.s., there exists the $g_{\Gamma} $-solution $y\in
\mathbf{D}_{\mathcal{F}}^{2}(0,T)$ with the terminal condition
$y_{T}=X $ (defined in Definition \ref{gGsol}). Moreover,
$g_{\Gamma}$-solution is the limit of a sequence of
$g^{n}$--solutions $y_{t}^{n}$ with $g^{n}=g+nd_{\Gamma }$, where
\begin{equation}
y_{t}^{n}=X+\int_{t}^{T}(g+nd_{\Gamma
})(s,y_{s}^{n},z_{s}^{n})ds-\int_{t}^{T}z_{s}^{n}dB_{s},  \label{gnGamma}
\end{equation}
with the convergence in the following sense:
\begin{eqnarray}
y_{t}^{n} &\nearrow &y_{t}\mbox{, with }\lim_{n\rightarrow \infty
}E[|y_{t}^{n}-y_{t}|^{2}]=0,\;\lim_{n\rightarrow \infty
}E\int_{0}^{T}|z_{t}-z_{t}^{n}|^{p}dt=0,\;  \label{approx1} \\
A_{t}^{n} &:&=\int_{0}^{t}(g+nd_{\Gamma
})(s,y_{s}^{n},z_{s}^{n})ds\rightarrow A_{t}\mbox{ weakly in }\mathbf{L}^{2}(\mathcal{%
F}_{t}),  \label{approx2}
\end{eqnarray}
where $z$ and $A$ are  corresponding martingale part and increasing
part of $y$, respectively.
\end{theorem}

\noindent\textbf{Proof. } By the comparison theorem of BSDE,
$y_t^n\leq y_t^{n+1}\leq y_t^{\prime }$. It follows that there
exists a $y\leq y^{\prime }$ such that, for each $t\in [0,T]$,
\[
y_t^1\leq y_t^n\nearrow y_t\leq y_t^{\prime }.
\]
Consequently, there exists a constant $C>0$, independent of $n$, such that
\[
E[\sup_{0\leq t\leq T}(y_t^n)^2]\leq C\,\,\,\, \hbox{and}\,\,\,\,
E[\sup_{0\leq t\leq T}(y_t^2)]\leq C.
\]
Thanks to the monotonic limit Theorem 2.1 in \cite{P99}, we can pass
limit on both sides of BSDE (\ref{gnGamma}) and obtain
\[
y_t=X+\int_t^Tg(s,y_s,z_s)ds+A_T-A_t-\int_t^Tz_sdB_s.
\]
On the other hand, by $E[(A_T^n)^2]=n^2E[(\int_0^td_{\Gamma
_s}(y_s^n,z_s^n)ds)^2]\leq C$, we have $ d_{\Gamma
_t}(y_t,z_t)\equiv0 $. \hfill$\Box$\medskip

\begin{remark}
\label{singular} From the approximation (\ref{gnGamma}) it is clear that, as
$n$ tends to $\infty ,$ the coefficient $g+nd_{\Gamma }$ tends to a singular
coefficient $g_{\Gamma }$ defined by
\[
g_{\Gamma }(t,y,z):=g(t,y,z)1_{\Gamma _{t}}(y,z)+\infty \times 1_{\Gamma
_{t}^{C}}(y,z).
\]
Thus, in the above theorem, the $g_{\Gamma }$-solution is also the
solution of BSDE with singular coefficient $g_{\Gamma }$.
\end{remark}

\begin{remark}
If the constraint $\Gamma $ is of the following form $\Gamma _{t}=(-\infty
,U_{t}]\times \mathbb{R}^{d}$, where $U_{t}\in \mathbf{L}^{2}(\mathcal{F}%
_{t})$, then the smallest $\Gamma $--constrained $g$--supersolution solution
with terminal condition $y_{T}=X$ exists, if and only if $d_{\Gamma
_{t}}(Y_{t},Z_{t})\equiv 0$, a.s. a.e., where $(Y,Z)$ is the solution of the
BSDE
\[
-dY_{t}=g(t,Y_{t},Z_{t})dt-Z_{t}dB_{t},\;\;t\in [0,T],\;\;Y_{T}=X.
\]
This follows easily by comparison theorem.
\end{remark}

We also have

\begin{theorem}[Comparison Theorem of $g_{\Gamma }$-solution]
\label{comp}We assume that $g^{1}$, $g^{2}$ satisfy (\ref{Lip}) and
$\Gamma ^{1}$, $\Gamma ^{2}$ satisfy (\ref{Gamma}). And suppose that
$\forall
(t,y,z)\in [0,T]\times \mathbb{R\times R}^{d}\mathbf{,}$%
\begin{equation}
X^{1}\leq X^{2},g^{1}(t,y,z)\leq g^{2}(t,y,z),\Gamma _{t}^{1}\supseteq
\Gamma _{t}^{2},  \label{con-compg}
\end{equation}
For $i=1,2$, Let $Y^{i}\in \mathbf{D}_{\mathcal{F}}^{2}(0,T)$ be the $%
g_{\Gamma ^{i}}^{i}$--solution with terminal condition $Y_{T}^{i}=X^{i}$.
Then we have
\[
Y_{t}^{1}\leq Y_{t}^{2},\mbox{ for }t\in [0,T],\;\mbox{a.s.}
\]
\end{theorem}

\noindent\textbf{Proof. } Consider the penalization equations for
the two constrained BSDE: for $n\in
\mathbf{N}$%
\begin{eqnarray}
y_{t}^{1,n}
&=&X^{1}+\int_{t}^{T}g^{1,n}(s,y_{s}^{1,n},z_{s}^{1,n})ds-%
\int_{t}^{T}z_{s}^{1,n}dB_{s},  \label{compp} \\
y_{t}^{2,n}
&=&X^{2}+\int_{t}^{T}g^{2,n}(s,y_{s}^{2,n},z_{s}^{2,n})ds-%
\int_{t}^{T}z_{s}^{2,n}dB_{s},  \nonumber
\end{eqnarray}
where
\begin{eqnarray*}
g^{1,n}(t,y,z) &=&g^{1}(t,y,z)+nd_{\Gamma _{t}^{1}}(y,z),\; \\
g^{2,n}(t,y,z) &=&g^{2}(t,y,z)+nd_{\Gamma _{t}^{2}}(y,z).
\end{eqnarray*}
From (\ref{con-compg}) we have $g^{1,n}(t,y,z)\leq g^{2,n}(t,y,z)$.
It follows from the classical comparison theorem of BSDE that
$y_{t}^{1,n}\leq y_{t}^{2,n}$. While as $n\rightarrow \infty $,
$y_{t}^{1,n}\nearrow y_{t}^{1} $ and $y_{t}^{2,n}\nearrow
y_{t}^{2}$, where $y^{1}$, $y^{2}$ are
the $g_{\Gamma }-$solutions of the BSDEs respectively. It follows that $%
y_{t}^{1}\leq y_{t}^{2}$, $0\leq t\leq T$.\hfill $\Box $\medskip

The comparison theorem is a powerful tool and useful concept in BSDE
Theorem (cf. \cite{EPQ}). Here let us recall the main theorem of
reflected BSDE and related comparison theorem for the case of lower
obstacle $L$. We do not repeat the case for the upper obstacle since
it is essentially the same. This result, obtained in \cite{PX2005},
is a generalized version of \cite{EKPPQ}, \cite{H} and \cite{LX} for
the part of existence, and \cite{HLM} for the part of comparison
theorem.

\begin{theorem}[Reflected BSDE and related Comparison Theory]
\label{CompR} We assume that the coefficient $g$ satisfies Lipschitz
condition (\ref{Lip}) and the lower obstacle $L$ satisfies
(\ref{barc}). Then, for each $X\in \mathbf{L}^{2}(\mathcal{F}_{T})$
with $X\geq L_{T}$ there exists
a unique triple $(y,z,A)\in \mathbf{D}_{\mathcal{F}}^{2}(0,T)\times \mathbf{L%
}_{\mathcal{F}}^{2}(0,T;\mathbb{R}^{d})\times (\mathbf{D}_{\mathcal{F}%
}^{2}(0,T))$, where $A$ is an increasing process, such that
\[
y_{t}=X+\int_{t}^{T}g(s,y_{s},z_{s})ds+A_{T}-A_{t}-\int_{t}^{T}z_{s}dB_{s}
\]
and the generalized Skorokhod reflecting condition is satisfied: for
each $L^{*}\in \mathbf{D}_{\mathcal{F}}^{2}(0,T)$ such that
$y_{t}\geq L_{t}^{*}\geq L_{t},$; $dP\times dt$ a.s., we have
\[
\int_{0}^{T}(y_{s-}-L_{s-}^{*})d\overline{A}_{s}=0,\;\mbox{a.s.,\ \
}
\]
Moreover, if a coefficient $g^{\prime }$ an obstacle $L^{\prime }$
and terminal condition $X^{\prime }$ satisfy the same condition as
$g$, $L$
and $X$ with for $\forall (t,y,z)\in [0,T]\times \mathbb{R\times R}^{d}%
\mathbf{,}$%
\[
X^{\prime }\leq X,g^{\prime }(t,y,z)\leq g(t,y,z),L_{t}^{\prime }\leq
L_{t},\,\,\,dP\times dt-\hbox{a.s.},
\]
and if the triple $(y^{\prime },z^{\prime },A^{\prime })$ is the
corresponding reflected solution, then we have
\[
Y_{t}^{\prime }\leq Y_{t},,\,\,\,\forall \,\,t\in
[0,T],\,\,\,\,\mbox{a.s.}
\]
and for each $0\leq s\leq t\leq T$,
\[
A_{t}^{\prime }\leq A_{t},\;\;A_{t}^{\prime }-A_{s}^{\prime }\leq
A_{t}-A_{s}.
\]
\end{theorem}

\textbf{Acknowledgment. }The first author thanks to Freddy Delbaen
for a fruitful discussion, after which we have understood an
interesting point of view of $g_\Gamma $--solution, noted in Remark
\ref{singular}.

\end{document}